\documentclass[journal]{IEEEtran}  

\IEEEoverridecommandlockouts                              

\usepackage{graphics} 
\usepackage{epsfig} 
\usepackage{mathptmx} 
\usepackage{times} 
\usepackage{amsmath} 
\usepackage{amssymb}  
\usepackage{soul}
\usepackage[mathcal]{euscript} 
\usepackage{mathtools}
\usepackage{mathrsfs}

\usepackage{graphicx}
\graphicspath{{Figures/}} 
\usepackage{caption} 
\usepackage{subfig} 

\usepackage{url}

\usepackage{xcolor}
\definecolor{USred}{rgb}{0.74,0.1,0.1}
\definecolor{USblue}{rgb}{0.2,0.2,0.7}
\definecolor{green1}{cmyk}{0.82,0,1,0.3}

\definecolor{Royalblue}{cmyk}{1,0.30,0.2,0.2}

\newcommand{\numberset}{\mathbb}
\newcommand{\NN}{\numberset{N}}

\newcommand{\RR}{\numberset{R}}
\newcommand{\CC}{\numberset{C}}
\newcommand{\TT}{\numberset{T}}
\newcommand{\ZZ}{\numberset{Z}}
\newcommand{\EE}{\mathbb{E}}

\newcommand{\B}{\mathcal{B}}
\newcommand{\C}{\mathcal{C}}

\newcommand{\Z}{\mathcal{Z}}

\newcommand{\V}{\mathcal{V}}
\newcommand{\G}{\mathcal{G}}
\newcommand{\I}{\mathcal{I}}
\newcommand{\J}{\mathcal{J}}

\newcommand{\N}{\mathcal{N}}
\newcommand{\Q}{\mathcal{Q}}
\newcommand{\K}{\mathcal{K}}

\newcommand{\argmin}{\operatornamewithlimits{argmin}}
\newcommand{\argmax}{\operatornamewithlimits{argmax}}
\renewcommand{\vec}{\boldsymbol}
\DeclareMathOperator{\tr}{tr}
\DeclareMathOperator{\rank}{rank}

\DeclareRobustCommand{\vect}[1]{
	\ifcat#1\relax
	\boldsymbol{#1}
	\else
	\mathbf{#1}
	\fi}

\newcommand\indep{\protect\mathpalette{\protect\independenT}{\perp}}
\def\independenT#1#2{\mathrel{\rlap{$#1#2$}\mkern2mu{#1#2}}}

\newtheorem{definition}{Definition}
\newtheorem{theorem}{Theorem}
\newtheorem{lemma}{Lemma}
\newtheorem{proposition}{Proposition}
\newtheorem{problem}{Problem}

\newtheorem{example}{Example}

\title{A Scalable Strategy for the Identification of Latent-variable Graphical Models}

\author{Daniele Alpago, Mattia Zorzi, Augusto Ferrante 
	\thanks{}
	\thanks{D. Alpago, M. Zorzi and A. Ferrante are with the Department of Information Engineering, University of Padova, Padova, Italy; email:	 
		{\tt\small alpagodani@dei.unipd.it} (D. Alpago)
		{\tt\small zorzimat@dei.unipd.it} (M. Zorzi)
		{\tt\small augusto@dei.unipd.it} (A. Ferrante)}%
	\thanks{}%
}

\begin{document}
	
\maketitle


\begin{abstract}
	In this paper we propose an identification method for latent-variable graphical models associated to autoregressive (AR) Gaussian stationary processes. The identification procedure exploits the approximation of AR processes through stationary reciprocal processes thus benefiting of the numerical advantages of dealing with block-circulant matrices. These advantages become more and more significant as the order of the process gets large. We show how the identification can be cast in a regularized convex program and we present numerical examples that compares the performances of the proposed method with the existing ones.
\end{abstract}


\section{INTRODUCTION}
\IEEEPARstart{T}{he} ideas behind graphical models have their origins in several scientific areas, such as statistical physics and genetics back at the beginning of the last century. However, only recent developments of such ideas allowed to employ graphical models in identification problems involving high dimensional data \cite{Lauritzen, Dahlh, LindqARMA, SongVan, SongConf, MLManSaid, MEManSaid, ChaParWill}. In this direction, particularly useful are \emph{sparse graphical models}, i.e. graphs with few edges that describe the interactions between a large number of variables. Such models have become very popular in the literature in the recent years because, beside giving a concise representation of the phenomenon under scruting, sparsity implies a limited number of model's parameters thus avoiding overfitting in the identification procedure.\\
Although the latter is a desirable property, enforcing sparsity in the identification procedure is not always the best choice, as it may prevent a sufficiently rich description of the underlying phenomenon. Indeed, in many practical situations, the presence of a few common, hidden behaviors between the variables of interest explaining the most part of the interactions between the observed variables can be crucial. The fact that a sparse graphical model is not able to describe the essential features of this kind of phenomena motivates the introduction of the so-called \emph{latent-variable graphical models}. The latter consist in a two-layer graph where the conditional dependence relations between the observed variables are mainly due to the latent variables (i.e. variables not accessible to observations): each latent-variable (on the top-layer) is connected to the majority of the observed variables (on the bottom-layer), making the latter a sparse subgraph. Since the number of latent-variables is small, the overall graph has a reduced number of edges. In the simplest possible setting, one can associate this kind of models to a Gaussian random vector \cite{ChaParWill}. The particular graphical structure translates in a sparse plus low-rank decomposition of its concentration matrix. In \cite{ChaParWill} the identification of the sparse and the low-rank part of the concentration matrix has been cast in a regularized maximum-likelihood optimization problem. A dynamic version of this problem, i.e. the identification of latent-variable graphical models for AR Gaussian processes, has been considered in \cite{ZorzSep} where the problem has been shown to be strictly connected to a maximum-entropy problem. As showed in \cite{ZorziADMM}, this identification problem can be effectively solved by an ADMM-type algorithm. The optimization procedure, however, involves the inversion and the eigenvalue decomposition of matrices whose dimension is proportional to the product of the order of the process by the dimension of the process, making the procedure numerically critical when the order of the AR process is high, as it happens, for example, when the AR process is an approximation of an ARMA one.\\
In this paper we consider the problem of identifying latent-variable graphical models for stationary Gaussian reciprocal processes. The latter are  periodic stationary processes \cite{CarFerrPav}, \cite{Rec1}, \cite{Rec2}, \cite{Rec3}, \cite{ringh2016multidimensional}, \cite{PICCI_LINDQUIST}, \cite{CFP-CDC05} and they have been proven to be a worthy approximation of Gaussian AR processes, provided that the  period $N$ is sufficiently large \cite{FerrLinAlg}, \cite{PICCI_LINDQUIST}. We will show that the proposed identification procedure is in fact an approximation of the maximum entropy and maximum likelihood identification paradigms proposed for the classical AR processes. The fact that stationary reciprocal process can be modeled by means of block-circulant matrices represents a big numerical advantage as the inversion and the eigenvalue decomposition of such matrices can be performed robustly \cite{ringh2015fast} making the proposed procedure attractive also for the identification AR processes of high order and hence for  ARMA processes.\\
The paper is organized as follows: In Section \ref{sec:notation} we fix the notation and we recall the fundamental results used in the rest of the paper. In Section \ref{sec:rp} we introduce reciprocal processes and we explain how they are related to AR processes. In Section \ref{sec:gm} we characterize graphical models associated to reciprocal processes while, in Section \ref{sec:id}, we propose a convex optimization problem for the identification of such models. Section \ref{sec:admm} is devoted to the ADMM formulation of the optimization problem and Section \ref{sec:ne} reports numerical experiments concerning the implementation of the proposed procedure. Finally, in Section \ref{sec:conc} we draw the conclusions.

\section{NOTATION and BACKGROUND}\label{sec:notation}
In this paper we will deal both with real matrices and with matrix-valued functions defined on the unit-circle $\TT:=\{e^{i\theta}:\,\theta\in[-\pi,\pi]\}$. For such functions we will omit the dependence on $\theta$ when it is clear from the context, i.e. we will write $F$ in place of $F(e^{i\theta})$. The rank of a matrix $G$ is denoted by $\rank(G)$ while the (normal) rank of any $\CC^{p\times p}$-valued analytic function $F$ defined on $\TT$, is defined as
\begin{equation}\label{eq:normrk}
\rank(F) := \max_{\theta\in[-\pi,\pi]}\rank(F(e^{i\theta})).
\end{equation}
In the same fashion, the following notations will be used indifferently in the case that $G$ is a $\CC^{p\times p}$-valued function defined on $\TT$ or a square constant matrix: $G^\top$ denotes the transpose of $G$, $G^*$ its transpose-conjugate and $\text{diag}(G)\in\CC^p$ denotes the vector whose entries are the diagonal elements of $G$. $\ker(G)$ indicates the kernel of $G$. $G>0$ and $G\ge 0$ denote that $G$ is a positive definite and, respectively, positive semidefinite. $\tr(G),\,\det(G)$ and $G^{-1}$ denote the trace of $G$, the determinant of $G$ and its inverse, respectively. 
$I_p$ denotes the identity matrix of order $p$. 

We define the cone
\[
\mathcal{S}_p := \{F\in\mathcal{H}_p: \Phi-\alpha\,I_p \ge0 \text{ a.e. on $\TT$, for some }\alpha>0\},
\]
where $\mathcal{H}_p$ is the space of square integrable coercive functions defined on the unit circle and taking values in the space of $p\times p$ Hermitian matrices. For any $F\in\mathcal{H}_p$ we will use equivalently the notations 
\[
\int_{-\pi}^\pi\,F(e^{i\theta}) \,\frac{d\theta}{2\pi},\quad\int\,F
\]
for the integral of $F$ over $[-\pi,\pi]$ with respect to the normalized Lebesgue measure on $\TT$. We define also the family of matrix pseudo-polynomials
\[
\Q_{p,n} :=\left\{ \sum_{k=-n}^n\,Q_k\,e^{i\theta k},\quad Q_{-k}=Q_k^\top\in\RR^{p\times p}  \right\}.
\]
For any sub-interval $(x_1,x_2):=\{x:\,x_1<x<x_2\}$ of an interval $(a,b)\subset\RR$, we denote with $(x_1,x_2)^c$ the complement set of $(x_1,x_2)$ in $(a,b)$. $\EE[\cdot]$ denotes the expectation operator.

In this paper we will always consider AR processes of order $n$ and reciprocal processes of period $N$, i.e. completely specified in a finite interval of length $N$. All such processes are understood with zero mean  throughout the paper. It will be always assumed that $N>2n$ and that $N$ is an even number. The case with $N$ odd can be dealt in a similar way. We define the vector space $\C\subset\RR^{mN\times mN}$ of the (real) symmetric, block-circulant matrices
\[
\vect{C} = \text{circ}\{C_0,C_1,\dots,C_{\frac{N}{2}-1},C_{\frac{N}{2}},C_{\frac{N}{2}-1}^\top,\dots,C_1^\top\},
\]
whose first block-column is composed by the $m\times m$ blocks $C_0,C_1,\dots,C_{\frac{N}{2}-1},C_\frac{N}{2},C_{\frac{N}{2}-1}^\top,\dots,C_1^\top$. The space $\C$ is endowed with the inner product $\left<\vect{C},\vect{D}\right>_\C:= \tr(\vect{C}^\top\vect{D})$. The \emph{symbol} of the block-circulant matrix $\vect{C}\in\C$ is defined as the $m\times m$ pseudo-polynomial
\begin{equation}\label{eq:symb}
\Phi(\zeta) := \sum_{k=0}^{N-1}\,C_k\,\zeta^{-k},\quad \text{ with } \quad C_k = C_{N-k}^\top\text{ for } k>\frac{N}{2},
\end{equation}\noindent
where $\zeta:=e^{i\frac{2\pi}{N}}$ is the $N$-th root of unity.
\begin{proposition}\label{prop:diagsymb}
	Let $\vect{C}$ be a block-circulant matrix with symbol $\Phi(\zeta)$ defined by \eqref{eq:symb}. Then
	\begin{equation}
	\vect{C} =
	\vect{F}^*
	\text{diag}\left\{\Phi(\zeta^0),\,\Phi(\zeta^1),\,\cdots,\,\Phi(\zeta^{N-1})\right\}
	\vect{F},
	\end{equation}
	where $\vect{F}$ is the (Fourier) unitary block-matrix
	{\small\begin{equation*}
		\vect{F}=
		\frac{1}{\sqrt{N}}
		\begin{bmatrix}
		\zeta^{-0 \cdot 0}I    & \zeta^{-0 \cdot 1}     & \cdots & \zeta^{-0 \cdot (N-1)}I    \\
		\zeta^{-1 \cdot 0}I    & \zeta^{-1 \cdot 1}I     & \cdots & \zeta^{-1 \cdot (N-1)}I    \\
		\vdots                   & \vdots                   & \ddots & \vdots                       \\
		\zeta^{-(N-1) \cdot 0}I & \zeta^{-(N-1) \cdot 1}I & \cdots & \zeta^{-(N-1) \cdot (N-1)}I \\
		\end{bmatrix}.
		\end{equation*}}
\end{proposition}
This is a classical result in the scalar case; technical details for the block-circulant case can be found, for instance, in \cite[page 6]{tesi}. We define the subspace $\B\subseteq\C$ of symmetric, banded block-circulant $mN\times mN$ matrices of bandwidth $n$, with $N>2n$, containing the matrices of the form
\begin{equation}\label{eq:bandmtx}
\vect{B} = \text{circ}\{B_0,B_1,\cdots,B_n,0,\cdots,0,B_n^\top,\cdots,B_1^\top\},
\end{equation}
that inherits the inner product defined on $\C$. Note that, according to definition \eqref{eq:symb}, the symbol of a banded matrix $\vect{B}\in\B$ is
\begin{equation*}
\Psi(\zeta) = \sum_{k=-n}^n\, B_k\,\zeta^{-k}, \qquad B_{-k} = B_k^\top.
\end{equation*}
The projection operator $\mathsf{P}_\mathscr{B}:\C\to\B$ is defined as 
\[
\mathsf{P}_\mathscr{B}(\vect{C}) :=\text{circ}\{C_0,C_1,\cdots,C_n,0,\cdots,0,C_n^\top,\cdots,C_1^\top\}.
\]
Given $\Omega=\{(i,j):\,i,j=1,\dots,m\}$, the projection operator $\mathsf{P}_\Omega:\C\to\C$ is defined such that $\mathsf{P}_\Omega(\vect{C})$ is a block-circulant matrix whose blocks have support $\Omega$.

\section{RECIPROCAL PROCESSES}\label{sec:rp}
Let $\{\vect{y}(k),\,k=1,2,\dots,N\}$, be an $m$-dimensional Gaussian stationary stochastic process defined on a finite interval $[1,N]$. For $k=1,\dots,N$, we have $\vect{y}(k) := [\vec{y}_1(k)\,\dots\,\vec{y}_m(k)]^\top\in\RR^m$, therefore the process is completely characterized by the random vector $\vect{y} := [\vec{y}_1(1)\,\dots\,\vec{y}_m(1)\,\dots\,\dots\,\vec{y}_1(N)\,\dots\,\vec{y}_m(N)]^\top\in\RR^{mN}$. In \cite{CarFerrPav} it has been shown that $\vect{y}$ is a restriction of a wide-sense stationary periodic process of period $N$ defined on the whole integer line $\ZZ$ if and only if the $mN\times mN$ covariance matrix $\vect{\Sigma}$ of $\vect{y}$ is symmetric block-circulant:
\begin{equation}\label{eq:circmtx}
\vect{\Sigma} =  \text{circ}\{\Sigma_0,\Sigma_1,\dots,\Sigma_\frac{N}{2},\dots,\Sigma_1^\top\},
\end{equation}
where $\EE[\vect{y}(i)\vect{y}(j)^\top] = \Sigma_{i-j}$, $i,j = 1,\dots,N$, are the covariance lags of the process such that $\Sigma_k = \Sigma_{N-k}^\top$ for $k>N/2$. In view of the above equivalence, we will denote with $\vect{y}$ both the wide-sense stationary periodic process defined in the whole line $\ZZ$ and its restriction, depending on the context. 
A particular class of stationary periodic processes is represented by reciprocal processes.

\begin{definition}
	$\vect{y}$ is a reciprocal process of order $n$ on $[1,N]$ if, for all $t_1,t_2\in[1,N]$, the random variables of the process in the interval $(t_1,t_2)\subset[1,N]$ are conditionally independent to the random variables in $(t_1,t_2)^c$, given the $2n$ boundary values $\vect{y}(t_1-n+1),\dots,\vect{y}(t_1),\vect{y}(t_2),\dots,\vect{y}(t_2+n-1)$, where the sums $t-k$ and $t+k$ are to be understood modulo $N$.
\end{definition}

The following result has been proved in \cite[Theorem 3.3]{CarFerrPav}: it states that a reciprocal process is completely specified by a block-circulant matrix whose inverse has a banded structure.

\begin{theorem}\label{thm:conc}
	A non-singular $mN \times mN$-dimensional matrix $\vect{\Sigma}$ is the covariance matrix of a periodic reciprocal process of order $n$ if and only if its inverse is a positive definite symmetric block-circulant matrix which is banded of bandwidth $n$, namely $\vect{\Sigma}^{-1}\in\B$.
\end{theorem}

\mbox{\\}

Let $\hat{\Sigma}_0,\dots,\hat{\Sigma}_n$ be given estimates of the first $n+1$ covariance lags $\Sigma_0,\dots,\Sigma_n$ of the underlying reciprocal process. In view of Theorem \ref{thm:conc}, the identification of a reciprocal process can be formulated as the following matrix completion problem.

\begin{problem}\label{pb:covext}
	Given the $n+1$ estimates $\hat{\Sigma}_0,\dots,\hat{\Sigma}_n$, compute a sequence $\Sigma_{n+1},\dots,\Sigma_\frac{N}{2}$, in such a way to form a symmetric, positive definite block-circulant matrix 
	\[
	\vect{\Sigma} = \text{circ}\{\hat{\Sigma}_0,\dots,\hat{\Sigma}_n,\Sigma_{n+1},\dots,\Sigma_\frac{N}{2},\dots,\Sigma_{n+1}^\top,\hat{\Sigma}_n^\top,\dots,\hat{\Sigma}_1^\top\},
	\]
	with $\Sigma^{-1}\in\B$.
\end{problem}
It has been shown in \cite{PICCI_LINDQUIST, CarFerrPav} that a particular solution to Problem \ref{pb:covext} is the one which solves the following maximum entropy problem:
\begin{equation}\label{op:mep}
\begin{aligned}
\argmax_{\vect{\Sigma}\in\C} &\quad \log\det\vect{\Sigma}\\
\text{subject to } &\quad \vect{\Sigma}>0\\
&\quad \mathsf{P}_{\mathscr{B}}(\vect{\Sigma}-\hat{\vect{\Sigma}}) = 0.
\end{aligned}
\end{equation}
whose dual problem has been proven to be
\begin{equation}\label{op:mepdu}
\begin{aligned}
\argmin_{\vect{X}\in\B} &\quad -\log\det\vect{X} + \left<\vect{X},\,\hat{\vect{\Sigma}}\right>_{\C}\\
\text{subject to } &\quad \vect{X}>0
\end{aligned}
\end{equation}
where $\hat{\vect{\Sigma}}\in\B$ is the symmetric, banded block-circulant matrix of bandwidth $n$,
\[
\hat \Sigma= \text{circ}\{\hat \Sigma_0,\hat \Sigma_1,\dots,\hat \Sigma_n,0,\dots,0,\hat \Sigma_n^\top,\dots,\hat \Sigma_1^\top\},
\]
containing the covariance lags estimated from the data and
the optimal value of dual variable $\vect{X}$ is indeed equal to $\Sigma^{-1}$, i.e. the inverse of the solution of \eqref{op:mep}. Strong duality between \eqref{op:mep} and \eqref{op:mepdu} implies that \eqref{op:mep} and \eqref{op:mepdu} are equivalent. In what follows we assume that $\hat{\vect{\Sigma}}>0$ as it is a necessary condition for Problem \eqref{op:mep} to be feasible. In the case that $\hat{\vect{\Sigma}}$ is not positive definite, we can consider a positive definite banded block-circulant matrix sufficiently close to $\hat{\vect{\Sigma}}$ which can be obtained by solving a structured covariance estimation problem, see \cite{structcov1}, \cite{structcov2}.

\subsubsection*{AR approximation}
Next we recall how reciprocal processes can be seen as an approximation of autoregressive (AR) processes. More precisely, let $\vec{\mathsf{y}} :=\{\vec{\mathsf{y}} (t):\,t\in\ZZ\}$ be an $m$-dimensional, AR, full-rank, Gaussian wide-sense stationary process of order $n$, 
\begin{equation}\label{eq:AR}
\sum_{k=0}^{n}\,B_k\,\vec{\mathsf{y}} (t-k) = \vect{e}(t),\qquad\vect{e}(t)\sim\N(0,I_m), \quad t\in\ZZ,
\end{equation}
and let $R_k := \EE[\vec{\mathsf{y}} (t)\vec{\mathsf{y}} (t-k)^\top]$, $k\in\ZZ$, be its $k$-th covariance lag. The spectrum of $\vec{\mathsf{y}} $ is the Fourier transform of the sequence $R_k$ with $k\in\ZZ$, i.e.
\begin{equation}\label{eq:psdAR}
\Phi(e^{i\theta}) = \sum_{k=-\infty}^{\infty}\,R_k\,e^{-i\theta k},\qquad R_{-k}=R_k^\top,\ \ \theta\in[-\pi,\pi].
\end{equation}
Suppose now that $T$ observations $\mathsf{y}(1),\dots,\mathsf{y}(T)$ of the process $\vec{\mathsf{y}}$ are available, and let
\begin{equation}\label{eq:estcov}
\hat{R}_k = \frac{1}{T}\sum_{t=k}^T\,\mathsf{y}(t)\mathsf{y}(t-k)^\top,\qquad k = 0,1,\dots,n,
\end{equation}
be estimates of the first $n+1$ covariance lags $R_0,\dots,R_n$. The identification of such a process can be cast to a \emph{covariance extension problem}.

\begin{problem}\label{pb:covextpb}
	Given $n+1$ estimates $\hat{R}_0,\,\dots,\,\hat{R}_n$, complete them with a sequence $R_{n+1},\,R_{n+2},\,\dots$ in such a way that the Fourier transform of the extended (infinite) sequence is a power spectral density.	
\end{problem}

A particular solution of Problem \ref{pb:covextpb} is the one proposed by J. P. Burg in \cite{BurgPhd}: choose $R_{n+1},\,R_{n+2},\,\dots$ maximizing the \emph{entropy rate} of the process, i.e. that solves the following optimization problem
\begin{equation}\label{eq:burg}
\begin{aligned}
\argmax_{\Phi\in\mathcal{S}_m} &\quad \int\,\log\det\Phi\\
\text{subject to } &\quad \int e^{i\theta k} \,\Phi = \hat{R}_k,\qquad k=0,1,\dots,n.
\end{aligned}
\end{equation}
The dual of \eqref{eq:burg} has been shown to be, see for instance \cite{ZORZI201587}:
\begin{equation}\label{eq:burgdual}
\begin{aligned}
\argmin_{\Phi^{-1}\in\Q_{m,n}} &\quad \int\,-\log\det\Phi^{-1} +\left<\Phi^{-1},\,\hat{\Phi}\right>\\
\text{subject to } &\quad \hat{\Phi}> 0
\end{aligned}
\end{equation}
where 
\begin{equation}\label{eq:correl}
\hat{\Phi}(e^{i\theta}) = \sum_{k=-n}^n\,\hat{R}_k\,e^{-i\theta k},\qquad \hat{R}_{-k} = \hat{R}_k^\top,
\end{equation}
is the truncated periodogram of the process $\vec{\mathsf{y}}$. These kind of problems have been extensively studied and generalized in the recent years, see for instance \cite{family_zorzi, Ringh2015MultidimensionalRC, KLapprox, FerMasPav,georgiou2006relative,byrnes2000new,BEL02}.\\

We recall that, for $N\to\infty$, Toeplitz matrices can be approximated arbitrarily well by circulant matrices \cite[Lemma 4.2]{Gray}; hence, for $N\to\infty$, Problem \ref{pb:covext} consists in searching a completion that leads to an infinite positive definite block-Toeplitz covariance matrix, i.e. such that the Fourier transform of the resulting extended sequence is a power spectral density. By Theorem 3.1 in \cite{FerrLinAlg}, for $N\to\infty$, Problem \eqref{op:mep} is the classical Burg's maximum entropy problem whose solution is an AR process of order $n$. In light of this observation, we can understand the reciprocal process associated to the solution of \eqref{op:mepdu} as an approximation of the AR process solution of the Burg's maximum entropy problem \eqref{eq:burg}. In the following sections we will exploit this approximation for the identification of latent-variable AR graphical models.\\

\begin{figure}[h!]\centering
	\includegraphics[scale=0.8]{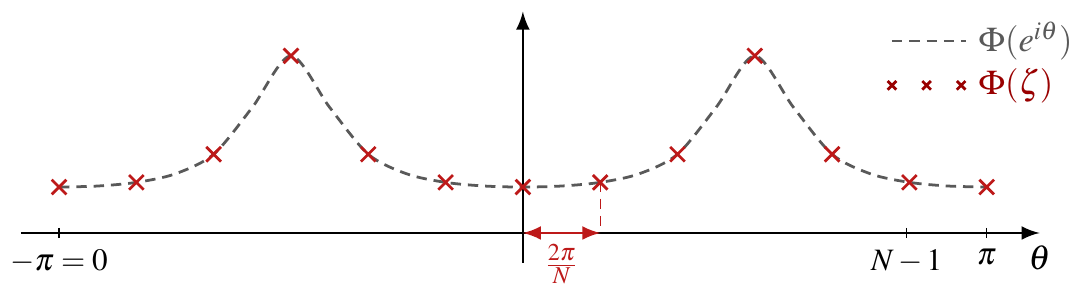}
	\caption{Spectrum $\Phi(e^{i\theta})$ and its sampled version $\Phi(\zeta)$ with $N=12$ samples.}\label{fig:recapx}
\end{figure}

\begin{figure}[h!]\centering
	\includegraphics[scale=0.75]{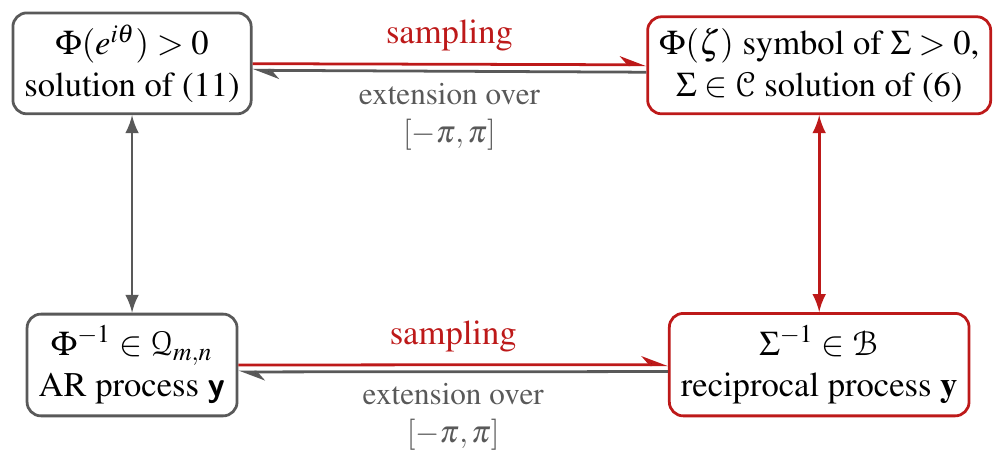}
	\caption{Schematic representation of the reciprocal approximation of the AR process in terms of solutions of problems \eqref{eq:burg} and \eqref{op:mep}.}\label{fig:recAR}
\end{figure}

The reciprocal approximation just explained has also an interesting interpretation in the frequency domain. Indeed, it corresponds to sampling the spectrum \eqref{eq:psdAR} of the AR process $\vec{\mathsf{y}} $, over the interval $[-\pi,\pi]$, with sample period $2\pi/N$, thus obtaining the \emph{symbol} of the covariance matrix of the corresponding reciprocal process:
\[
\Phi(\zeta) = \sum_{k=0}^{N-1}\, \Sigma_k \,\zeta^{-k},\qquad \Sigma_{k}=\Sigma_{N-k}^\top\text{ for } k>\frac{N}{2}.
\]
Figure \ref{fig:recapx} illustrates this relation. According to Proposition \ref{prop:diagsymb}, the covariance matrix $\vect{\Sigma}$ of the reciprocal process $\vect{y}$ that approximates $\vec{\mathsf{y}} $ writes as
\begin{equation}\label{eq:covrec}
\vect{\Sigma} = \vect{F}^*\text{circ}\{ \Phi(\zeta^0),\,\Phi(\zeta^1),\dots,\Phi(\zeta^{N-1}) \}\vect{F},
\end{equation}
hence, its inverse
\begin{equation}\label{eq:invcovrec}
\vect{\Sigma}^{-1} = \vect{F}^*\text{circ}\{ \Phi(\zeta^0)^{-1},\,\Phi(\zeta^1)^{-1},\dots,\Phi(\zeta^{N-1})^{-1} \}\vect{F},
\end{equation}
can be robustly computed by inverting the $N$ blocks $\Phi(\zeta^0),\,\Phi(\zeta^1),\dots,\Phi(\zeta^{N-1})$, all of size $m\times m$. As a final remark, we recall that eigevalues and eigenvectors of circulant matrices can be robustly computed as well, thanks to the availability of closed-form formulas, see for instance \cite{Gray}.%
\vspace{2mm}
\\
As highlighted by the frequency-domain interpretation, the goodness of the approximation  depends on the regularity of the spectrum: the larger is the rate of variation of the spectrum, the larger $N$ has to be chosen in order to get a good approximation of the AR process. The frequency-domain interpretation makes even more explicit the relationship between Burg's maximum entropy problem \eqref{eq:burg} and Problem \eqref{op:mep}: provided that the number of samples $N$ is sufficiently large, by sampling the spectrum solution of \eqref{eq:burg} we obtain an approximation of the matrix $\vect{\Sigma}$ solution of \eqref{op:mep}; viceversa, the symbol of $\vect{\Sigma}$ can be extended over the whole interval $[-\pi,\pi]$ in order to approximate the  solution of \eqref{eq:burg}. Figure \ref{fig:recAR} summarizes this bi-directional relationships.

\section{GRAPHICAL MODELS}\label{sec:gm}
Consider a Gaussian random vector $\vect{x}\sim\N(0,\Sigma)$ taking values in $\RR^m$, where $\Sigma=\Sigma^\top>0$ so that the concentration matrix $K=[k_{ij}]:=\Sigma^{-1}$ is well-defined. If we denote the components of $\vect{x}$ as $\vec{x}_1,\dots,\vec{x}_m$, for any $i\ne j$, we have that $x_i$ is \emph{conditionally independent} from $x_j$ given the remaining random variables $x_k$, $k\ne i,j$, i.e.
\begin{equation}\label{eq:cishn}
x_i \indep x_j \,|\, \{x_k\}_{k\ne i,j},
\end{equation}
if and only if the element $k_{ij}$ in position $(i,j)$ of the concentration matrix $K$ is equal to zero. Formally,
\begin{equation}\label{eq:ci}
   x_i \indep x_j \,|\, \{x_k\}_{k\ne i,j}\quad\iff\quad k_{ij}=0.
\end{equation}
The previous relation allows to construct an undirected graph $\G=(V,E)$, with $V=\{1,\dots,m\}$ and $E\subset V\times V$, associated to the random vector $\vect{x}$ by taking the components $\vec{x}_1,\dots,\vec{x}_m$ of $\vect{x}$ as nodes and such that the edges reflect the conditional dependence relations between the random variables, i.e.
\begin{equation}\label{eq:edg}
   (i,j)\notin E \quad\iff\quad x_i \indep x_j \,|\, \{x_k\}_{k\ne i,j}.
\end{equation}
The graph $\G$ is called the \emph{graphical model} associated to $\vect{x}$ and it gives a visual representation of the conditional dependence relations between the components of $\vect{x}$. Observe that $\G$ is completely characterized by the sparsity pattern of the concentration matrix of the random vector. 

A characterization of conditional independence can be given also in the dynamic setting.
In particular, we consider an $m$-dimensional, Gaussian, wide-sense stationary AR process $\vect{x}$ described by a model like \eqref{eq:AR}. For any index set $I\subset V$, define 
\begin{equation*}
\mathcal{X}_I := \text{span}\{\vec{x}_j(t):\,j\in I,\,t\in\ZZ\},
\end{equation*}
as the closure of the set containing all the finite linear combinations of the variables $\vec{x}_j(t)$. For any $i\ne j$, the notation
\begin{equation*}
\mathcal{X}_{\{i\}} \indep \mathcal{X}_{\{j\}} \,|\,\mathcal{X}_{V\setminus\{i,j\}}
\end{equation*}
generalizes \eqref{eq:cishn} and it means that for all $t_1,t_2$,  $\vec{x}_i(t_1)$ and $\vec{x}_j(t_2)$ are conditionally independent given the space linearly generated by $\{\vec{x}_k(t),\,k\in V\setminus\{i,j\}, t\in\ZZ\}$. One can prove that
\begin{equation}\label{eq:cidin}
\mathcal{X}_{\{i\}} \indep \mathcal{X}_{\{j\}} \,|\,\mathcal{X}_{V\setminus\{i,j\}}\quad\iff\quad [\Phi(e^{i\theta})^{-1}]_{ij}=0,
\end{equation}
for any $\theta\in[-\pi,\pi]$, see \cite{Dahlh,LindqARMA}, which is the natural generalization of \eqref{eq:ci}. Accordingly, we can construct the undirected graph $\G=(V,E)$ representing the conditional dependence relations between the components of the process $\vect{x}$ by defining the set of edges as follows:
\begin{equation}
(i,j)\notin E \quad\iff\quad \mathcal{X}_{\{i\}} \indep \mathcal{X}_{\{j\}} \,|\,\mathcal{X}_{V\setminus\{i,j\}}.
\end{equation}
In this framework, the graph $\G$ is completely characterized by the sparsity pattern of the inverse power spectral density of the process. Identification of sparse graphical models of reciprocal processes have been studied in \cite{nostroELETTERS}.\\
In many practical situations there is the presence of a few common, latent, behaviors between the variables of interest that are responsible of the most part of the interactions between the observed variables and that cannot be captured by considering only a sparse model structure. This leads to a particular type of graphical models called \emph{latent-variable graphical models} or \emph{sparse plus low-rank graphical models}, \cite{ChaParWill}. Such models admit a two-layer graphical structure in which the nodes in the upper layer stand for the (few) latent-variables, while the nodes in the bottom layer represent the observed variables.\\
Latent-variable graphical models associated to Gaussian random vectors have been considered in \cite{ChaParWill} and then generalized in \cite{ZorzSep} to  AR stochastic processes. The latter, say $\vec{\mathsf{z}}:=\{\vec{\mathsf{z}}(t),\,t\in\ZZ\}$, is assumed to be of the form $\vec{\mathsf{z}} = [\vec{\mathsf{y}}^\top\,\vec{\mathsf{x}}^\top]^\top$ where $\vec{\mathsf{y}}$ is the $\RR^m$-valued process containing the observed variables while $\vec{\mathsf{x}}$ is the process containing $l$ latent variables. Let $\Phi_{\vec{\mathsf{y}}}$ denotes the spectral density of $\vec{\mathsf{y}}$. Under the assumptions that $l\ll m$ and the dependence relations among the observed variables are mostly through the latent variables, we have the decomposition
\begin{equation}\label{eq:S+Lpsd}
   \Phi_{\vec{\mathsf{y}}}^{-1} = \Gamma - \Lambda,
\end{equation}
where $\Gamma>0$ is sparse and its support reflects the conditional dependencies among the observed variables, while $\Lambda\ge0$ is low-rank and its rank   equals  the number $l$ of latent-variables.

We are now ready to extend the previous results for Gaussian reciprocal processes. Let $\vect{z}:=[\vect{y}^\top\,\vect{x}^\top]^\top$ be a Gaussian, periodic, reciprocal process of order $n$ defined on the interval $[1,N]$, where $\vect{y}$ plays the role of the $m$-dimensional observed process and $\vect{x}$ is the $l$-dimensional latent process, respectively. The covariance matrix $\vect{\Sigma}_\vect{z}$ of $\vect{z}$ and its inverse can be partitioned as
\begin{equation}\label{eq:ShurSig}
\renewcommand{\arraystretch}{1.3}
	\vect{\Sigma}_\vect{z} = \left[
	\begin{array}{c|c}
	\vect{\Sigma}_\vect{y} & \vect{\Sigma}_{\vect{y}\vect{x}}\\\hline
	\vect{\Sigma}_{\vect{y}\vect{x}}^\top & \vect{\Sigma}_\vect{x}
	\end{array}\right],
	\qquad
	\vect{\Sigma}_\vect{z}^{-1} = \left[
	\begin{array}{c|c}
	\vect{S} & \vect{A}\\\hline
	\vect{A}^\top & \vect{R}
	\end{array}\right],
\end{equation}
where $\vect{\Sigma}_\vect{y}\in\C$ and $\vect{\Sigma}_\vect{x}\in\C_l$ are the covariance matrices of $\vect{y}$ and $\vect{x}$, respectively. Here, $\C_l$ denotes the vector space of block-circulant, symmetric matrices as $\C$, except that the blocks have dimension $l\times l$.
Applying the Schur complement, we obtain the relation
\begin{equation}\label{eq:S+L}
\vect{\Sigma}_\vect{y}^{-1} = \vect{S}-\vect{L},
\end{equation}
where $\vect{S}>0$ is the concentration matrix of process $\vect{y}$ conditioned on $\vect{x}$, and $\vect{L}\ge 0$ is defined as $\vect{L}:=\vect{A}\,\vect{R}^{-1}\,\vect{A}^\top$. In order to ensure that, according to Theorem \ref{thm:conc}, $\vect{\Sigma}_\vect{y}^{-1}\in\B$ we assume both $\vect{S}$ and $\vect{L}$ to be symmetric, block-circulant, banded of bandwidth $n$, i.e.
\begin{equation}\label{eq:SLstruct}
\begin{aligned}
\vect{S} &=\text{circ}\{S_0,S_1,\dots,S_n,0,\dots,0,S_n^\top,\dots,S_1^\top\},\\
\vect{L} &=\text{circ}\{L_0,L_1,\dots,L_n,0,\dots,0,L_n^\top,\dots,L_1^\top\}.
\end{aligned}
\end{equation}
By construction, the matrix $\vect{L}$ has rank equal to the number of latent variables $l$, therefore under the assumption that $l\ll m$, it is a low-rank matrix. If $\vect{S}$ is a sparse matrix, then we will refer to \eqref{eq:S+L} as \emph{sparse plus low-rank decomposition} of $\vect{\Sigma}_\vect{y}^{-1}$ which is the analogue of \eqref{eq:S+Lpsd} for reciprocal processes. It remains to show that an appropriate sparsity pattern of $\vect{S}$ reflects that the dependence relations among observed variables are mostly through the few latent variables. For this purpose, let $\vect{y}_i:=[\vec{y}_i(1) \dots \vec{y}_i(N)]^\top$, $i=1,\dots,m$, be the $i$-th component of the process $\vect{y}$ and let $\vect{x}_j:=[\vec{x}_j(1) \dots \vec{x}_j(N)]^\top$, $j=1,\dots,l$, be the $j$-th component of the process $\vect{x}$. 
Although the components of the reciprocal processes are defined \emph{for any} $k\in\ZZ$, by periodicity it is sufficient to impose conditional independence only for $k\in[1,N]$. We assume that the blocks $S_0, S_1,\dots,S_n$ of $\vect{S}$ have common support $\Omega\subseteq\{(i,j):\,i,j=1,\,\dots,\,m\}$ namely,
\begin{equation}\label{eq:cs}
(S_k)_{ij} = (S_k)_{ji} = 0,\qquad k=0,\,\dots\,,n,\quad\forall\,(i,j)\in\Omega^c,
\end{equation}
where $\Omega$ is the set of pairs that contains all the $(i,\,i),\,i=1,\dots,m$. By property \eqref{eq:ci}, equation \eqref{eq:cs} is equivalent to
\begin{equation}\label{eq:cs1}
\begin{aligned}
\vec{y}_i(t_1)\indep\vec{y}_j(t_2)\mid\,&\{{y}_h(s),\,h\ne i,j,\,\,s=1,\dots,N,\\
& \vec{y}_i(s_1),\,s_1\ne t_1,\,\vec{y}_j(s_2),\,s_2\ne t_2,\,\vect{x}\},
\end{aligned}
\end{equation}
for any $t_1,\,t_2\in[1,N]$ and for any pair $(i,j)\in\Omega^c$.
\begin{proposition}\label{prop:tesi}
	Condition \eqref{eq:cs1} is equivalent to
	\begin{equation}\label{eq:condindp}
	\vec{y}_i(t_1)\indep\vec{y}_j(t_2)\mid 
	\{\vec{y}_h(s),\,h\ne i,j,\,\,s=1,\dots,N,\,\vect{x}\}
	\end{equation}
	for any $t_1,\,t_2\in[1,N]$ and for any $(i,j)\in\Omega^c$.
	\begin{IEEEproof}
    The proof exploits basic results of the theory of Hilbert spaces of second-order random variables, see for instance \cite[Chapter 2]{lindpicLSS}. First of all, let
    \[
    \vect{\epsilon} := \begin{bmatrix}
    \vec{\epsilon}_i\\\vec{\epsilon}_j
    \end{bmatrix}=
    \begin{bmatrix}
    \vect{y}_i\\\vect{y}_j
    \end{bmatrix} - 
    \EE\left[ \begin{bmatrix}
    \vect{y}_i\\\vect{y}_j
    \end{bmatrix}\,\biggm\vert\,\vect{y}_h(s),\,h\ne i,j,\,s=1,\dots,N,\vect{x}\right]
    \]
    denotes the error affecting the projection of $[\vect{y}_i^\top\,\,\vect{y}_j^\top]^\top$ onto the subspace generated by $\{\vect{y}_h(s),\,h\ne i,j,\,s=1,\dots,N,\vect{x}\}$, for any $t_1,\,t_2\in[1,N]$ and for any $(i,j)\in\Omega^c$.  It can be shown that $\vect{\epsilon}$ is a zero-mean, Gaussian, random vector. Accordingly, proving \eqref{eq:condindp} is equivalent to prove that
    \begin{equation}\label{eq:epsorth}
    \EE\left[\vec{\epsilon}_i\,\vec{\epsilon}_j^\top\right] = 0\qquad\iff\qquad
    \vec{\epsilon}_i(t_1)\,\indep\,\vec{\epsilon}_j(t_2)
    \end{equation}
    for any $t_1,\,t_2\in[1,N]$ and for any $(i,j)\in\Omega^c$, \cite{lindpicLSS}. Let now $\Pi$ be a permutation matrix that permutes the rows of $\vect{z}=[\vect{y}^\top\,\vect{x}^\top]^\top$ in order to obtain
    \[
    \bar{\vect{z}} := \Pi\,\vect{z} =
    \left[
    \begin{array}{c}
    \vect{y}_i\\
    \vect{y}_j\\\hline
    \vect{y}_{h\ne i,j}\\
    \vect{x}
    \end{array}
    \right]
    =\left[
    \begin{array}{c}
    \bar{\vect{z}}_1\\\hline
    \bar{\vect{z}}_2
    \end{array}\right],
    \]
    where $\vect{y}_{h\ne i,j}$ is the vector containing the random variables $\vect{y}_h(s),\,h\ne i,j,\,s=1,\dots,N$. We partition the covariance matrix $\vect{\Sigma}_{\bar{\vect{z}}}$ of $\bar{\vect{z}}$ as
    \[\renewcommand{\arraystretch}{1.3}
    \vect{\Sigma}_{\bar{\vect{z}}} = \left[
    \begin{array}{c|c}
    \vect{\Sigma}_{\bar{\vect{z}}_1} & \vect{\Sigma}_{\bar{\vect{z}}_1\bar{\vect{z}}_2}\\\hline
    \vect{\Sigma}_{\bar{\vect{z}}_2\bar{\vect{z}}_1} & \vect{\Sigma}_{\bar{\vect{z}}_2}
    \end{array}\right],
    \] 
    where $\vect{\Sigma}_{\bar{\vect{z}}_1}$ and $\vect{\Sigma}_{\bar{\vect{z}}_2}$ are the covariance matrices of $\bar{\vect{z}}_1$ and $\bar{\vect{z}}_2$, respectively. It is well known that its inverse can be partitioned conformably as
    \[\renewcommand{\arraystretch}{1.3}
    \vect{\Sigma}_{\bar{\vect{z}}}^{-1} = \Pi\,\Sigma_\vect{z}^{-1}\,\Pi^\top 
    =
    \left[
    \begin{array}{c|c}
    \bar{\vect{S}}& *\\\hline
    * & *
    \end{array}\right],
    \]
    where
    \begin{equation}
    	\bar{\vect{S}}:=\left(\vect{\Sigma}_{\bar{\vect{z}}_1} - \vect{\Sigma}_{\bar{\vect{z}}_1\bar{\vect{z}}_2}\vect{\Sigma}_{\bar{\vect{z}}_2}^{-1}\vect{\Sigma}_{\bar{\vect{z}}_2\bar{\vect{z}}_1}\right)^{-1}
    \end{equation}   
    is a permuted version of matrix $\vect{S}$, according to the permutation matrix $\Pi$. By construction, the Schur complement formula applied on $\vect{\Sigma}_{\bar{\vect{z}}}$ gives 
    \begin{equation}\label{eq:seps}
    \Sigma_{\vect{\epsilon}} = \vect{\Sigma}_{\bar{\vect{z}}_1} - \vect{\Sigma}_{\bar{\vect{z}}_1\bar{\vect{z}}_2}\vect{\Sigma}_{\bar{\vect{z}}_2}^{-1}\vect{\Sigma}_{\bar{\vect{z}}_2\bar{\vect{z}}_1} =\bar{\vect{S}}^{-1},
    \end{equation}
    that relates the covariance matrix $\Sigma_{\vect{\epsilon}}$ of the projection error $\vect{\epsilon}$ to the covariance matrix $\vect{\Sigma}_{\bar{\vect{z}}}$ of $\bar{\vect{z}}$. Condition \eqref{eq:cs} is equivalent to say that $\bar{\vect{S}}$, and therefore $\bar{\vect{S}}^{-1}$, is block-diagonal. Accordingly, by \eqref{eq:seps}, $\vect{\Sigma}_{\vec{\epsilon}}$ is block-diagonal, i.e. $\vec{\epsilon}_i$ and $\vec{\epsilon}_j$ are independent, which is equivalent to \eqref{eq:epsorth} as we wanted to prove.
	\end{IEEEproof}
\end{proposition}

The above result reflects the fact that the random variables $\{\vec{y}_i(s_1),\vec{y}_j(s_2),\,s_1\ne t_1,\,s_2\ne t_2\}$ do not play any role in the conditioning \eqref{eq:cs1}. Moreover, the group-sparsity condition \eqref{eq:cs} translates in the fact that the conditional dependence relations between the observed variables are mainly due to the few latent variables. Accordingly, condition \eqref{eq:condindp} represents the reciprocal counterpart of condition \eqref{eq:ci}. We conclude that $\vect{\Sigma}_{\vect{z}}^{-1}$ in \eqref{eq:ShurSig} together with \eqref{eq:cs} define an undirected graph for the Gaussian random vector $\vect{z}$ which admits a two-layer structure where
\begin{itemize}
	\item[-] The nodes in the upper-layer represent the $l$ variables of the latent-process $\vect{x}_1,\dots,\vect{x}_l$ while the nodes in the bottom-layer represent the $m$ variables of the observed process $\vect{y}_1,\dots,\vect{y}_m$.
	\item[-] The edges are given by the entries of the concentration matrix $\vect{\Sigma}_\vect{z}^{-1}$. In particular, the edge $(i,j)$, between two vectors $\vect{y}_i$ and $\vect{y}_j$, $i\ne j$, is described by
	\[
	\left[ (S_0)_{ij}\,\,(S_1)_{ij}\,\,\dots\,\,(S_n)_{ij}\,\, 0 \,\, \dots \,\, 0 \,\, (S_n)_{ji}\,\,\dots\,\, (S_1)_{ji}\right].
	\]
\end{itemize}

\begin{example}
	Consider the case in which $N=2$, $m=7$, $l=2$, and suppose that the graphical model associated to the vector $\vect{z}$ is the one depicted in Figure \ref{fig:2layergraph}.
	\begin{figure}[h!]\centering
		\includegraphics[scale=0.9]{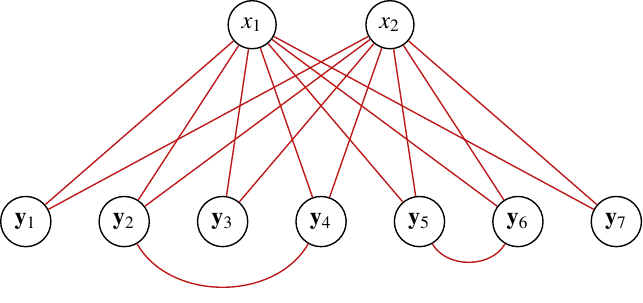}
		\caption{Example of a latent-variable graphical model: $\vec{x}_1,\,\vec{x}_2$ are the latent-variables and $\vect{y}_1,\vect{y}_2,\dots,\vect{y}_7$ are the manifest variables.}\label{fig:2layergraph}
	\end{figure}\\
	In this case, the concentration matrix of vector $\vect{z}$ will have the structure \eqref{eq:ShurSig} with
	\begin{equation*}\renewcommand{\arraystretch}{1.3}
	\vect{S}=\left[
	\begin{array}{c|c}
	S_0 & S_1^\top\\\hline
	S_1 & S_0
	\end{array}\right],\quad
	S_0,\,S_1\in\RR^{7\times7},
	\end{equation*}
	and $R$ is a $2\times 2$ matrix. The presence of an edge between $\vect{y}_2$ and $\vect{y}_4$ implies that at least one of the two elements $(S_0)_{24}$ and $(S_1)_{24}$ is different from zero. Similar arguments holds for the edge between $\vect{y}_5$ and $\vect{y}_6$. Thus, $\Omega=\{ (i,i):\,i=1,\dots,7\}\cup\{(2,4),\,(5,6)\}$.
\end{example}

\section{IDENTIFICATION of LATENT-VARIABLE RECIPROCAL GRAPHICAL MODELS}\label{sec:id}
The problem of identifying a latent-variable graphical model associated to a Gaussian random vector has been firstly considered in \cite{ChaParWill} where the solution is obtained by solving a regularized maximum likelihood problem. In \cite{ZorzSep} the problem has been extended to a dynamic setting, by considering an AR Gaussian process. More precisely, in \cite{ZorzSep} a regularized version of Problem \eqref{eq:burgdual} that relies on the sparse plus low-rank decomposition of the inverse of the observed spectrum in \eqref{eq:S+Lpsd}, has been considered:
\begin{equation}\label{eq:SpLcl}
\begin{aligned}
\argmin_{\Gamma,\,\Lambda\in\Q_{m,n}} &\quad \int\,-\log\det(\Gamma - \Lambda) + \langle\Gamma - \Lambda,\,\hat{\Phi}_{\vec{\mathsf{y}}}\rangle\\
&\quad +\gamma_S\,\phi_1(\Gamma)+\gamma_L\,\phi_*(\Lambda)\\
\text{subject to } &\quad \Gamma - \Lambda> 0\\
&\quad \Lambda \ge 0.
\end{aligned}
\end{equation}
Here, $\gamma_S,\,\gamma_L>0$ are the regularization parameters that balance the effects of the two regularizers $\phi_1$ and $\phi_*$ inducing sparsity and low-rank on $\Gamma$ and $\Lambda$, respectively, while $\hat{\Phi}_{\vec{\mathsf{y}}}$ is the truncated periodogram of the observed process $\vec{\mathsf{y}}$. In this section we propose a procedure for the identification of a latent-variable graphical model associated to an AR Gaussian process that exploits the approximation of an AR process through a reciprocal process in the sense explained in Section \ref{sec:rp}. Recalling that a latent-variable graphical model of a reciprocal process is characterized by \eqref{eq:S+L}, the system identification problem can be stated as follows.
\begin{problem}\label{pb:SpLid}
	Consider an $m$-dimensional AR process $\vec{\mathsf{y}}$ and let $\hat{R}_0,\dots,\hat{R}_n$ be the estimates of the first $n+1$ covariance lags of $\vec{\mathsf{y}}$ computed as in \eqref{eq:estcov}. Set $\Sigma_0:=\hat{R}_0,\dots,\Sigma_n:=\hat{R}_n$. Compute the blocks $\Sigma_{n+1},\dots,\Sigma_\frac{N}{2}$ of the block-circulant covariance matrix 
	$\vect{\Sigma}_\vect{y}=\text{circ}\{\Sigma_0,\Sigma_1,\dots,\Sigma_{\frac{N}{2}-1},\Sigma_{\frac{N}{2}},\Sigma_{\frac{N}{2}-1}^\top,\dots,\Sigma_1^\top\}$ such that $\vect{\Sigma}_\vect{y}^{-1} = \vect{S} - \vect{L}$, 
	where $\vect{S}>0$ and $\vect{L}\ge 0$ are as in \eqref{eq:SLstruct} with $S_0,\dots,S_n$ having the smallest possible common support $\Omega$, as in \eqref{eq:cs}, and the $\rank$ of $\vect{L}$ is as small as possible.
\end{problem}

We stress the fact that only samples of the observed processes are available. Clearly, the matrix $\vect{\Sigma}_\vect{y}$ solving Problem \ref{pb:SpLid} is the covariance of the reciprocal  process $\vect{y}$ approximating the observed process $\vec{\mathsf{y}}$. Since we are going to identify a model for a reciprocal process, we can exploit the maximum entropy dual problem \eqref{op:mepdu} recalled in Section \ref{sec:rp}. It is worth noting that the support $\Omega$ is not known in advance, thus it has to be estimated from the data. In order to do that, inspired by \cite{SongVan}, we consider the following regularizer
\[
h_\infty(\vect{S}) = \sum_{k>h}\,\max\left\{ |(S_0)_{hk}|,2\max_{j=1,\dots,n}|(S_j)_{hk}|,2\max_{j=1,\dots,n}|(S_j)_{kh}|\right\}.
\]
The latter is a generalization of the $\ell^\infty$-norm used to induce group-sparsity on vectors, and it is used to enforce on $\vect{S}$ the group sparsity in \eqref{eq:cs}. The trace (as a tractable proxy of the nuclear norm) is used instead for inducing low-rankness in $\vect{L}$. Therefore, the paradigm for the estimation of the sparse plus low-rank decomposition of the concentration matrix $\vect{\Sigma}_\vect{y}^{-1}$ now directly follows from \eqref{op:mepdu} by setting $\vect{X}=\vect{S}-\vect{L}$, with $\vect{L} \ge0$, and by adding the regularizers just introduced:
\begin{equation}\label{eq:primal}
\begin{aligned} 
\argmin_{\vect{S},\vect{L}\in\B} &\quad-\log\det(\vect{S}-\vect{L}) + \left<\hat{\vect{\Sigma}}_{\vec{\mathsf{y}}},\,\vect{S}-\vect{L}\right>_\C\\
& \quad + \lambda_S\,h_\infty(\vect{S})+\lambda_L\tr(\vect{L})\\
\text{subject to } &\quad \vect{S}-\vect{L}>0,\quad\vect{L}\ge 0
\end{aligned}
\end{equation}
where $\lambda_L,\,\lambda_S>0$ are the two regularization parameters and 
\[
\hat{\vect{\Sigma}}_{\vec{\mathsf{y}}}= \text{circ}\{\hat R_0,\hat R_1,\dots,\hat R_n,0,\dots,0,\hat R_n^\top,\dots,\hat R_1^\top\}
\]
is the symmetric, banded block-circulant matrix of bandwidth $n$, containing the covariance lags estimated from the observations. As a further motivation, observe that Problem \eqref{eq:primal} is precisely the reciprocal counterpart of Problem \eqref{eq:SpLcl} considered in \cite{ZorzSep}. By replacing $\vect{S}$ with $\vect{X}:=\vect{S}-\vect{L}$, it becomes
\begin{equation}\label{eq:primal1}
\begin{aligned} 
\argmin_{\vect{X},\vect{L}\in\B} &\quad-\log\det(\vect{X}) + \tr(\hat{\vect{\Sigma}}_{\vec{\mathsf{y}}}\,\vect{X}) + \lambda_S\,h_\infty(\vect{X}+\vect{L})+\lambda_L\tr(\vect{L})\\
\text{subject to } &\quad \vect{X}>0,\quad\vect{L}\ge 0.
\end{aligned}
\end{equation}
We address the previous constrained optimization problem using the Lagrange multipliers theory. In doing that we add a new dummy variable $\vect{Y}$
\begin{equation}\label{eq:primal2}
\begin{aligned} 
\argmin_{\substack{\vect{X}\in\C\\\vspace{1mm}\vect{Y},\vect{L}\in\B}} &\quad-\log\det(\vect{X}) + \tr(\hat{\vect{\Sigma}}_{\vec{\mathsf{y}}}\,\vect{X}) + \lambda_S\,h_\infty(\vect{Y})+\lambda_L\tr(\vect{L})\\
\text{subject to } &\quad \vect{X}>0,\quad\vect{L}\ge 0\\
&\quad \vect{Y} = \vect{X} + \vect{L}.
\end{aligned}
\end{equation}
The Lagrangian function for this problem is 
\begin{equation}\label{eq:lag}  
\begin{aligned}
\mathcal{L}(\vect{X},\vect{Y},\vect{L},\vect{V},\vect{Z}) = &-\log\det(\vect{X}) + \left<\hat{\vect{\Sigma}}_{\vec{\mathsf{y}}},\,\vect{X}\right>_{\C} + \lambda_S\,h_\infty(\vect{Y})\\
&+\lambda_L\tr(\vect{L}) - \left<\vect{V},\,\vect{L}\right>_{\C} + \left<\vect{Z},\vect{X}+\vect{L}-\vect{Y}\right>_{\C}
\end{aligned}
\end{equation}
where, $\vect{V}\in\B$, because $\vect{L}\in\B$, and $\vect{V}\ge 0$, while $\vect{Z}\in\C$. After simple computations we have
\begin{align*}
\mathcal{L}(\vect{X},\vect{Y},\vect{L},\vect{V},\vect{Z}) = &-\log\det(\vect{X}) + \left<\hat{\vect{\Sigma}}_{\vec{\mathsf{y}}}+\vect{Z},\,\vect{X}\right>_{\C}\\
&+\left<\lambda_L  I_{mN} - \vect{V}+\vect{Z},\,\vect{L}\right>_{\C}\\
&+\lambda_S\,h_\infty(\vect{Y})- \left<\vect{Z},\,\vect{Y}\right>_{\C}.
\end{align*}
The dual objective function is the infimum over $\vect{X},\,\vect{Y}$ and $\vect{L}$ of the Lagrangian. The unique term on $\mathcal{L}$ that depends on $\vect{Y}$ is
$ \lambda_S\,h_\infty(\vect{Y}) - \left<\vect{Z},\,\vect{Y}\right>_{\C}$. The latter
is bounded below if and only if
\begin{align}
&\text{diag}(Z_j)=0,\quad  j=0,\,\dots,\,n, \label{eq:minY1}\\
&2|(Z_0)_{kh}|+\sum_{j=1}^n\,|(Z_j)_{kh}|+|(Z_j)_{hk}|\le\frac{\lambda_S}{N},\quad k>h,\label{eq:minY2}
\end{align}
in which case the infimum is zero. Accordingly,
\[
\inf_\vect{Y}\,\mathcal{L}=
\left\{
\begin{split}
&-\log\det(\vect{X}) + \left<\hat{\vect{\Sigma}}_{\vec{\mathsf{y}}}+\vect{Z},\,\vect{X}\right>_{\C} +\left<\lambda_L  I_{mN} - \vect{V}+\vect{Z},\,\vect{L}\right>_{\C}\\ &\text{if}\quad\eqref{eq:minY1},\,\eqref{eq:minY2}\text{ hold,}\\\\
&-\infty\quad\text{ otherwise.}
\end{split}
\right.
\]
The only term that depends on $\vect{L}$ is $\left<\lambda_L  I_{mN} - \vect{V}+\vect{Z},\,\vect{L}\right>_{\C}$. Recalling that $\vect{L}, \vect{V}\in\B$, by using the linearity of the projection operator $\mathsf{P}_\mathscr{B}$, we have that
\begin{equation}\label{eq:Vconstr}
\left<\lambda_L I_{mN} - \vect{V}+\vect{Z},\,\vect{L}\right>_{\C} =\left<\lambda_L I_{mN} - \vect{V}+\mathsf{P}_\mathscr{B}(\vect{Z}),\,\vect{L}\right>_{\C}
\end{equation}
which is linear in $\vect{L}$, and therefore it is bounded below if and only if 
\begin{equation} \label{eq:minL}
\lambda_L  I_{mN} - \vect{V}+\mathsf{P}_\mathscr{B}(\vect{Z})=0.
\end{equation}
In this case, the minimum of \eqref{eq:Vconstr} is zero. Accordingly,
\[
\inf_{\vect{Y},\vect{L}}\,\mathcal{L}=
\left\{
\begin{split}
&-\log\det(\vect{X}) + \left<\hat{\vect{\Sigma}}_{\vec{\mathsf{y}}}+\vect{Z},\,\vect{X}\right>_{\C}\\
&\text{if} \quad \eqref{eq:minY1},\,\eqref{eq:minY2},\,\eqref{eq:minL}\text{ hold,}\\\\
&-\infty\quad\text{ otherwise.}
\end{split}
\right.
\]
If $\eqref{eq:minY1},\,\eqref{eq:minY2},\,\eqref{eq:minL}$ hold, it remains to minimize the strictly convex function
\[
\bar{\mathcal{L}}(\vect{X}) := \inf_{\vect{Y},\vect{L}}\,\mathcal{L} = -\log\det(\vect{X}) + \left<\hat{\vect{\Sigma}}_{\vec{\mathsf{y}}}+\vect{Z},\,\vect{X}\right>_{\C}
\]
over the cone of the symmetric, positive definite, banded block-circulant matrices. Observe that, for any $\vect{Z}\in\C$, any $\hat{\vect{\Sigma}}_{\vec{\mathsf{y}}}\in\B$, and for any sequence $\vect{X}_k>0$ converging to a singular matrix, 
\[
\lim_{k\rightarrow \infty}\bar{\mathcal{L}}(\vect{X}_k)=\infty.
\]
Accordingly, we can assume that the solution lies in the interior of the cone so that a necessary and sufficient condition for $\vect{X}_o$ to be a minimum point for $\bar{\mathcal{L}}$ is that its first Gateaux derivative computed at $\vect{X}=\vect{X}_o$ is equal to zero in every direction $\vect{\delta}\vect{X}$, namely
\begin{equation}\label{eq:ddcond}
\vect{\delta}\bar{\mathcal{L}}(\vect{X}_o;\vect{\delta}\vect{X}) = \tr\left[ \left(-\vect{X}_o^{-1} + \hat{\vect{\Sigma}}_{\vec{\mathsf{y}}} + \vect{Z}\right)\vect{\delta}\vect{X}\right] = 0,
\qquad \forall\,\vect{\delta}\vect{X}\in\C.
\end{equation}
Notice that $\bar{\mathcal{L}}$ is bounded below if and only if
\begin{equation}\label{eq:pdcstr}
\hat{\vect{\Sigma}}_{\vec{\mathsf{y}}}+\vect{Z} > 0,
\end{equation}
therefore condition \eqref{eq:ddcond} is satisfied if and only if $\vect{X}_o = (\hat{\vect{\Sigma}}_{\vec{\mathsf{y}}}+\vect{Z})^{-1}$. Hence,
\[
\inf_{\vect{Y},\vect{L},\vect{X}}\,\mathcal{L}=
\left\{
\begin{split}
&\log\det(\hat{\vect{\Sigma}}_{\vec{\mathsf{y}}}+\vect{Z}) + mN,\\
&\text{ if } \eqref{eq:minY1},\,\eqref{eq:minY2},\,\eqref{eq:minL},\,\eqref{eq:pdcstr}\text{ hold,}\\\\
&-\infty\qquad\text{ otherwise.}
\end{split}
\right.
\]
Therefore, the dual problem of \eqref{eq:primal} is
\begin{equation}
\begin{aligned} 
\argmin_{\vect{V}\in\B,\,\vect{Z}\in{\C}} &\quad -\log\det(\hat{\vect{\Sigma}}_{\vec{\mathsf{y}}}+\vect{Z}) - mN\\
\text{subject to } & \quad V\geq 0,\,\eqref{eq:minY1},\,\eqref{eq:minY2},\,\eqref{eq:minL},\,\eqref{eq:pdcstr}.
\end{aligned}
\end{equation}
Notice that we can remove the variable $\vect{V}$. Indeed, recalling that $\vect{V}\ge 0$, the constraint \eqref{eq:minL} becomes $\lambda_LI_{mN} +\mathsf{P}_\mathscr{B}(\vect{Z}) = \vect{V} \ge 0$. Accordingly, the dual problem takes the form
\begin{equation}\label{eq:dual}
\begin{aligned} 
\argmin_{\vect{Z}\in{\C}} &\quad -\log\det(\hat{\vect{\Sigma}}_{\vec{\mathsf{y}}} +\vect{Z}) - mN\\
\text{subject to } &\quad \eqref{eq:minY1},\,\eqref{eq:minY2},\,\eqref{eq:pdcstr}\\
&\quad \lambda_LI_{mN} + \mathsf{P}_\mathscr{B}(\vect{Z}) \ge 0.\\
\end{aligned}
\end{equation}
\begin{proposition}
	Under the assumption that $\hat{\vect{\Sigma}}_{\vec{\mathsf{y}}}\in\B$ and $\hat{\vect{\Sigma}}_{\vec{\mathsf{y}}}>0$, Problem \eqref{eq:dual} admits a unique solution.
	\begin{IEEEproof}
		Define $f(\vect{Z}):=\log\det(\hat{\vect{\Sigma}}_{\vec{\mathsf{y}}}+\vect{Z})$. Let
		\[
		  \Q :=\left\{ \vect{Z}\in\C\,|\,\eqref{eq:minY1},\,\eqref{eq:minY2},\,\eqref{eq:pdcstr}\text{ and }\lambda_LI_{mN} + \mathsf{P}_\mathscr{B}(\vect{Z}) \ge 0\text{ hold}\right\}
		\]
		be the set of constraints of Problem \eqref{eq:dual}. First of all, notice that constraints \eqref{eq:minY1} and \eqref{eq:minY2} ensure that $\Q$ is a bounded subset of ${\C}$. Indeed, the entries of any $\vect{Z}\in\Q$ are bounded by $\lambda_S/N$ so that $\|\vect{Z}\|_{\C}<\infty$ for any $\vect{Z}\in\Q$. Let now $(\vect{Z}^{(k)})_{k\in\NN}$ be a generic sequence of elements of $\Q$ converging to some $\bar{\vect{Z}}\in{\C}$, such that $\hat{\vect{\Sigma}}_{\vec{\mathsf{y}}}+\bar{\vect{Z}}\ge0$ is singular. Then
		\[
		\lim_{k\to\infty}\,-\log\det(\hat{\vect{\Sigma}}_{\vec{\mathsf{y}}}+\vect{Z}^{(k)})=+\infty,
		\]
		and therefore $\vect{Z}^{(k)}$ is not an infimizing sequence. Hence, we can restrict the research of the minimum to the closed subset of $\Q$ defined by 
		\begin{align*}
		\bar{\Q}:=\{ \vect{Z}\in\C\,|\,&\hat{\vect{\Sigma}}_{\vec{\mathsf{y}}}+\vect{Z}\ge\epsilon I_{mN},\,\eqref{eq:minY1},\,\eqref{eq:minY2}\\
		&\text{and } \lambda_LI_{mN} + \mathsf{P}_\mathscr{B}(\vect{Z}) \ge 0 \text{ hold}\}
		\end{align*}
		with $\epsilon>0$ small enough. By what we have shown till now, the function $f$ is continuous on the compact set $\bar{\Q}$ and therefore it admits at least one minimum point. Since $f$ is strictly convex, the minimum is also unique.
	\end{IEEEproof}
\end{proposition}

\begin{proposition}
	Under the assumption that $\hat{\vect{\Sigma}}_{\vec{\mathsf{y}}}\in\B$ and $\hat{\vect{\Sigma}}_{\vec{\mathsf{y}}}>0$, Problem \eqref{eq:primal1} admits a solution $(\vect{X}_o,\vect{L}_o)$ and $\vect{X}_o$ is unique.
\end{proposition}
\begin{IEEEproof}
	Notice that Problem \eqref{eq:primal1} is a strictly feasible convex optimization problem (for instance, pick $\vect{X}=I_{mN}$ and $\vect{L}=0$). Accordingly, Slater's condition holds, hence strong duality holds between \eqref{eq:primal1} and its dual. The strong duality between problems \eqref{eq:primal1} and \eqref{eq:dual} and the existence of a unique optimum $\vect{Z}_o$ for the dual problem \eqref{eq:dual}, imply that there exists a unique $\vect{X}_o\in\B$ so that
	$
	\vect{X}_o = \left(\hat{\vect{\Sigma}}_{\vec{\mathsf{y}}} + \vect{Z}_o\right)^{-1}
	$
	which solves the primal problem \eqref{eq:primal1}.\\	
	It remains to show that there exists an $\vect{L}_o\in\B$ that solves the optimization problem 
	\begin{equation}\label{eq:optpbL}
	\begin{aligned} 
	\argmin_{\vect{L}\in\B} &\quad \lambda_S\,h_\infty(\vect{X}_o+\vect{L})+\lambda_L\tr(\vect{L})\\
	\text{subject to } &\quad \vect{L}\ge 0.
	\end{aligned}
	\end{equation}
	Notice that, the objective function in \eqref{eq:optpbL} is continuous. Since $\vect{L}=0$ is a feasible point, the problem is equivalent to find $\vect{L}\in\B$ that minimizes $\lambda_S\,h_\infty(\vect{X}_o+\vect{L})+\lambda_L\tr(\vect{L})$ over the set
	\[
	\K:= \left\{ \vect{L}\in\B\,\bigg|\, \vect{L}\ge0,\,\,\lambda_S\,h_\infty(\vect{X}_o+\vect{L})+\lambda_L\tr(\vect{L})\le\lambda_S\,h_\infty(\vect{X}_o)\right\}.
	\]
	It is easy to see that $\K$ is a closed and bounded (and thus compact) subset of $\B$. Hence, by Weierstrass' Theorem, Problem \eqref{eq:optpbL} admits a solution $\vect{L}_o$. At this point we can conclude that the primal problem \eqref{eq:primal1} admits a solution $(\vect{X}_o,\,\vect{L}_o)$.
\end{IEEEproof}
%

\subsection{Interpretations}
In the remaining of this section we will show how Problem \eqref{eq:primal} can be interpreted either as a regularized maximum-likelihood problem or as a dual of a maximum entropy problem.
\subsubsection*{Maximum likelihood interpretation} The reciprocal approximation of AR processes illustrated in Section \ref{sec:rp} allows to interpret Problem \eqref{eq:primal} as a regularized (conditional) maximum likelihood problem. Indeed, in the following we will show that the fitting function in \eqref{eq:primal}, i.e.
\begin{equation}\label{eq:primobjnoreg}
-\log\det(\vect{S}-\vect{L}) + \tr\left(\hat{\vect{\Sigma}}_{\vec{\mathsf{y}}}\,(\vect{S}-\vect{L})\right),
\end{equation}
is the approximation of the (conditional) negative log-likelihood of the AR process \eqref{eq:AR} that should be understood in the sense explained in Section \ref{sec:rp}. Following \cite{songsiri2010graphical}, consider the observed AR process $\vec{\mathsf{y}}$ whose spectrum is denoted by  $\Phi_{\vec{\mathsf{y}}}$,
and suppose that $T$ observations $\mathsf{y}(1),\dots,\mathsf{y}(T)$ of the process are available. The conditional likelihood of the process $\vec{\mathsf{y}}$ is defined as the likelihood function associated to the conditional distribution of $\mathsf{y}(n+1),\mathsf{y}(n+2),\dots,\mathsf{y}(n+T)$ given  $\mathsf{y}(1),\dots,\mathsf{y}(n)$. Let  
\[
\vect{T}_n := \text{Toepl}\{\hat{R}_0,\,\hat{R}_1\,,\cdots,\,\hat{R}_n\}
\]
be the block-Toeplitz matrix having in the first rows the estimates of the first $n+1$ covariance lags of the process $\hat{R}_0,\hat{R}_1,\dots,\hat{R}_n$ computed as in \eqref{eq:estcov}. For $T$ large enough, the conditional negative log-likelihood function of the AR process can be well approximated by
\begin{equation*}
	\ell(B) := -(T-n)\,\log\det B_0 + \frac{T-n}{2}\,\tr(B\,\vect{T}_n\,B^\top)
\end{equation*}
where $B:=[B_0\,\,B_1\,\,\cdots\,\,B_n]$ is the $(n+1)m$-dimensional vector containing the coefficients of the process. Applying Jensen's formula, it turns out that
\begin{equation*}
	\log\det B_0 = \frac{1}{2}\int \log\det\Phi_{\vec{\mathsf{y}}}(e^{i\theta}),
\end{equation*}
moreover, if $\hat{\Phi}_{\vec{\mathsf{y}}}(e^{i\theta})$
is the truncated periodogram of the AR process in \eqref{eq:correl}, it is easy to see that
\begin{equation*}
\int \hat{\Phi}_{\vec{\mathsf{y}}}(e^{i\theta})\,e^{i\theta k} = \hat{R}_{-k} = \hat{R}_k^\top. 
\end{equation*}
Accordingly, the approximated conditional negative log-likelihood can be rewritten as
\begin{equation}\label{eq:clpsd}
	\ell(B) = \frac{T-n}{2}\int \log\det\Phi_{\vec{\mathsf{y}}}(e^{i\theta}) + \tr\left[\hat{\Phi}_{\vec{\mathsf{y}}}(e^{i\theta})\,\Phi_{\vec{\mathsf{y}}}(e^{i\theta})^{-1}\right].
\end{equation}
A natural way to approximate \eqref{eq:clpsd} is to approximate the integral with a finite sum, i.e. to discretize the interval $[-\pi,\pi]$. This is precisely the frequency interpretation of the reciprocal approximation explained in Section \ref{sec:rp} that consists in sampling the spectrum of the process to obtain the corresponding symbol (see Figure \ref{fig:recapx}). In fact, considering as sample frequency $\Delta\theta=2\pi/N$, the Backward Euler approximation leads to the discrete approximation
\begin{equation*}
   \ell(B) \simeq\frac{T-n}{2}\frac{\Delta\theta}{2\pi}\,\sum_{k=0}^{N-1}\log\det\Phi_{\vec{\mathsf{y}}}(e^{i\theta_k}) + \tr\left[\hat{\Phi}_{\vec{\mathsf{y}}}(e^{i\theta_k})\Phi_{\vec{\mathsf{y}}}(e^{i\theta_k})^{-1}\right].
\end{equation*}
where $\theta_k = k\,\Delta\theta-\pi$. The conditional log-likelihood can now be rewritten straightforward in terms of symbols as
\begin{equation*}
	\ell(B) \simeq \frac{T-n}{2N}\left[\sum_{k=0}^{N-1}\log\det\Phi_{\vec{\mathsf{y}}}(\zeta^k) + \tr\sum_{k=0}^{N-1}\hat{\Phi}_{\vec{\mathsf{y}}} (\zeta^k)\Phi_{\vec{\mathsf{y}}} (\zeta^k)^{-1}\right].
 \end{equation*}
Observe now that $\Phi_{\vec{\mathsf{y}}} (\zeta)$ is precisely the symbol of the block-circulant covariance matrix $\vect{\Sigma}_\vect{y}$ of the reciprocal process $\vect{y}$ approximating the process $\vec{\mathsf{y}}$ and $\hat{\Phi}_{\vec{\mathsf{y}}} (\zeta)$ is the symbol of the block-circulant matrix $\hat{\vect{\Sigma}}_{\vec{\mathsf{y}}}$ in Problem \eqref{eq:primal}. Accordingly, form Proposition \ref{prop:diagsymb}, it follows that
\begin{equation*}
  \ell(B) \simeq \frac{T-n}{2N}\left[ -\log\det\vect{\Sigma}_\vect{y}^{-1}+\tr\left( \hat{\vect{\Sigma}}_{\vec{\mathsf{y}}}\,\vect{\Sigma}_\vect{y}^{-1}\right) \right].
\end{equation*}
Since $\vect{\Sigma}_\vect{y}^{-1}=\vect{S}-\vect{L}$, this is precisely (up to a scaling factor) equal to \eqref{eq:primobjnoreg}.

\subsubsection*{Maximum entropy interpretation} We will show that Problem \eqref{eq:primal} can be interpreted a regularized version of the dual of a maximum entropy problem, see \cite{P-F-SIAM-REV} for a general overview of these problems. Consider the regularized solution $(\vect{S}_o,\,\vect{L}_o)$ of \eqref{eq:primal} and let $\Omega$ be the support of $\vect{S}_o$, i.e. $\vect{S}_o$ satisfies \eqref{eq:cs}. Since $\vect{L}_o\in\B$ is so that $\vect{L}_o\ge0$ and $\rank\vect{L}_o=lN\ll mN$, there exists
\begin{equation*}
	\vect{G} = \text{circ}\{G_0,G_1,\dots,G_n,0,\dots,0\}
\end{equation*}
such that $G_k\in\RR^{m\times l}$ and $\vect{L}_o=\vect{G}^\top\vect{G}$. Accordingly, we can consider a modified version of Problem \eqref{eq:primal} where the regularizers are replaced by the corresponding hard-constraints $\vect{S}\in\V_\Omega$ and $\vect{L}\in\V_G$, where $\V_\Omega :=\{\vect{S}\in\C:\,\mathsf{P}_{\Omega^c}(\vect{S})=0\}$ and $\V_G :=\{ \vect{G}^\top(I_N\otimes H)\vect{G}:\,H\in\RR^{l\times l},\,H=H^\top\}$ is such that $\V_G\subseteq\B$. Thus, the resulting problem is
\begin{equation}\label{eq:hardprimal}
\begin{aligned} 
\argmin_{\vect{S},\vect{L}\in\B} &\quad-\log\det(\vect{S}-\vect{L}) + \left<\hat{\vect{\Sigma}}_{\vec{\mathsf{y}}},\,\vect{S}-\vect{L}\right>_\C\\
\text{subject to } &\quad \vect{S}-\vect{L}>0,\,\,\vect{L}\ge 0,\\
&\quad \vect{S}\in\V_\Omega,\,\,\vect{L}\in\V_G.
\end{aligned}
\end{equation}

\begin{proposition}
	The primal of Problem \eqref{eq:hardprimal} is
	\begin{equation}\label{eq:pdME}
	\begin{aligned} 
	\argmax_{\vect{\Sigma}_\vect{y}\in\C} &\quad \log\det\vect{\Sigma}_\vect{y}\\
	\text{subject to } &\quad \mathsf{P}_\Omega\mathsf{P}_\B(\vect{\Sigma}_\vect{y}-\hat{\vect{\Sigma}}_{\vec{\mathsf{y}}})=0,\\
	&\quad \vect{E}^*\,\vect{G}(\vect{\Sigma}_\vect{y}-\hat{\vect{\Sigma}}_{\vec{\mathsf{y}}})\vect{G}^\top\vect{E}\ge 0,
	\end{aligned}
	\end{equation}
	where $\vect{E}^*:=\frac{1}{\sqrt{N}}\left[I_l\quad 0 \quad \cdots \quad 0\right]$.
	\begin{IEEEproof}
		We derive the dual of Problem \eqref{eq:pdME}. Observing that $\vect{E}=\vect{F}^*\vect{1}$ where $\vect{1}:=\frac{1}{\sqrt{N}}\left[I_l\quad I_l \quad \cdots \quad I_l\right]^\top$, the Lagrangian of Problem \eqref{eq:pdME} writes as
		\begin{align*}
		\mathcal{L}(\vect{\Sigma}_\vect{y},\vect{W},H) = &\log\det{\vect{\Sigma}_\vect{y}} + \left<\mathsf{P}_{\Omega\cup\B}(\hat{\vect{\Sigma}}_{\vec{\mathsf{y}}}-\vect{\Sigma}_\vect{y}),\,\vect{W}\right>_\C\\
		&+\left<\vect{1}^\top\vect{F}\,\vect{G}(\vect{\Sigma}_\vect{y}-\hat{\vect{\Sigma}}_{\vec{\mathsf{y}}})\vect{G}^\top\vect{F}^*\vect{1},\,H\right>_\C,
		\end{align*}
		where $\vect{W}\in\C$, $H\in\RR^{l\times l}$ is a positive semidefinite symmetric matrix, and $\mathsf{P}_{\Omega\cup\B}(\vect{S})=\mathsf{P}_\Omega\mathsf{P}_\B(\vect{S})$. The last term of the Lagrangian can be rewritten as
		\begin{align*}
		&\tr\left[\vect{F}(\vect{\Sigma}_\vect{y}-\hat{\vect{\Sigma}}_{\vec{\mathsf{y}}})\vect{F}^*\,\,\vect{F}\vect{G}^\top\vect{F}^*\vect{1}H\vect{1}^\top\,\,\vect{F}\vect{G}\vect{F}^*\right]\\
		&= \tr\left[\vect{F}(\vect{\Sigma}_\vect{y}-\hat{\vect{\Sigma}}_{\vec{\mathsf{y}}})\vect{F}^*\,\,\vect{F}\vect{G}^\top\vect{F}^*(I_N\otimes H)\,\,\vect{F}\vect{G}\vect{F}^*\right]\\
		&= \tr\left[(\vect{\Sigma}_\vect{y}-\hat{\vect{\Sigma}}_{\vec{\mathsf{y}}})\,\,\vect{G}^\top(I_N\otimes H)\,\vect{G}\right],
		\end{align*}
		where we have exploited the fact that $\vect{F}(\vect{\Sigma}_\vect{y}-\hat{\vect{\Sigma}}_{\vec{\mathsf{y}}})\vect{F}^*$ and $\vect{F}\vect{G}\vect{F}^*$ are block-diagonal matrices and the fact that $\vect{F}^*(I_N\otimes H)\vect{F}=I_N\otimes H$. Accordingly,
		\begin{align*}
		\mathcal{L}(\vect{\Sigma}_\vect{y},\vect{W},H) = &\log\det{\vect{\Sigma}_\vect{y}} + \left<\hat{\vect{\Sigma}}_{\vec{\mathsf{y}}}-\vect{\Sigma}_\vect{y},\,\mathsf{P}_{\Omega\cup\B}(\vect{W})\right>_\C\\
		&+\left<\vect{\Sigma}_\vect{y}-\hat{\vect{\Sigma}}_{\vec{\mathsf{y}}},\,\,\vect{G}^\top(I_N\otimes H)\,\vect{G}\right>_\C\\
		=&\log\det{\vect{\Sigma}_\vect{y}} + \left<\hat{\vect{\Sigma}}_{\vec{\mathsf{y}}}-\vect{\Sigma}_\vect{y},\,\vect{S}-\vect{L}\right>_\C,
		\end{align*}
		where $\vect{S}:=\mathsf{P}_{\Omega\cup\B}(\vect{W})$ belongs to $\V_\Omega$ and $\vect{L}:=\vect{G}^\top(I_N\otimes H)\,\vect{G}\ge0$ belongs to  $\V_G\subseteq\B$, i.e. they satisfy all the constraints in \eqref{eq:hardprimal}.
		Similar arguments as the ones used to prove formula \eqref{eq:ddcond}, allow us to assert that a necessary and sufficient condition for $\vect{\Sigma}_o$ to be a minimum point for $\mathcal{L}$ is that its first Gateaux derivative computed at $\vect{\Sigma}_\vect{y}=\vect{\Sigma}_o$ is equal to zero in every direction $\vect{\delta}\vect{\Sigma}$, namely
		\begin{equation*}
		\vect{\delta}\mathcal{L}(\vect{\Sigma}_o;\vect{\delta}\vect{\Sigma}) = \tr\left[ \left(\vect{\Sigma}_o^{-1} -\vect{S} + \vect{L}\right)\vect{\delta}\vect{\Sigma}\right] = 0,
		\qquad \forall\,\vect{\delta}\vect{\Sigma}\in\C.
		\end{equation*}
		By assumption we have that $\vect{S}-\vect{L}>0$ thus, the substitution of the optimum $\vect{\Sigma}_o=(\vect{S}-\vect{L})^{-1}$ in the Lagrangian $\mathcal{L}$ leads precisely to the objective function in \eqref{eq:hardprimal}.
		\end{IEEEproof}
\end{proposition}
Some observations on the two constraints of \eqref{eq:pdME} are in order. The first constraint $\mathsf{P}_\Omega\mathsf{P}_\B(\vect{\Sigma}_\vect{y}-\hat{\vect{\Sigma}}_{\vec{\mathsf{y}}})=0$ fixes the entries corresponding to the indexes in $\Omega$ of the first $n+1$ lags of the reciprocal process. Concerning the second constraint, let $\Psi(\zeta)$ and $\Phi_\vect{y}(\zeta)$ be the symbols of $\vect{G}$ and $\vect{\Sigma}_\vect{y}$, respectively. By Proposition \ref{prop:diagsymb}, we have that
\begin{equation}\label{eq:filter}
\vect{E}^*\,\vect{G}\,\vect{\Sigma}_\vect{y}\,\vect{G}^\top\vect{E}=\frac{1}{N}
\sum_{k=0}^{N-1}\,\Psi(\zeta^k)\,\Phi_\vect{y}(\zeta^k)\,\Psi(\zeta^k)^*,
\end{equation}
which is the covariance of the output of the $m\times l$ filter $   \Psi(\zeta)=\sum_{k=0}^n\,G_k\,\zeta^{-k}$ fed with the reciprocal process $\vect{y}$. Accordingly, the second constraint in \eqref{eq:pdME} states that the covariance matrix of the process at the output of the filter is lower-bounded by $\vect{E}^*\,\vect{G}\,\hat{\vect{\Sigma}}_{\vec{\mathsf{y}}}\,\vect{G}^\top\vect{E}$. We conclude that Problem \eqref{eq:pdME} can be seen as the reciprocal counterpart of the maximum entropy problem \cite{ZorzSep},
\begin{equation}\label{eq:merelax}
\begin{aligned}
\argmax_{\Phi_{\vec{\mathsf{y}}}\in\mathcal{S}_m} &\quad \int\,\log\det\Phi_{\vec{\mathsf{y}}}\\
\text{subject to } &\quad \left(\int e^{i\theta k} \,\Phi_{\vec{\mathsf{y}}}-\hat{R}_k\right)_{pq}=0,\quad \substack{k=0,1,\dots,n\\(p,q)\in\Omega}\\
&\quad \int \Psi\,(\Phi_{\vec{\mathsf{y}}} - \hat{\Phi}_{\vec{\mathsf{y}}})\,\Psi^*\ge 0,
\end{aligned}
\end{equation}
where $\Psi(e^{i\theta})=\sum_{k=0}^n\,G_k\,e^{-i\theta k}$. Indeed, the second constraint in \eqref{eq:merelax} can be approximated with the backward Euler approximation with sample frequency $\Delta\theta=2\pi/N$ obtaining \eqref{eq:filter}.

\section{ALTERNATING DIRECTION METHOD of MULTIPLIERS} \label{sec:admm}
The solution of Problem \eqref{eq:dual} requires the joint enforcement of the constraints \eqref{eq:minY1}, \eqref{eq:minY2} and $\lambda_L\,I_{mN}+\mathsf{P}_\B(\vect{Z})\geq 0$, which may be a difficult task. In this section we will use the alternating direction methods of multipliers (ADMM) \cite{BoydADMM} to solve Problem \eqref{eq:dual} by showing that the constraints can be separated and each one can be enforced in an alternating way.\\
First of all observe that, by defining the variable $\vect{P}:=\lambda_LI_{mN} + \mathsf{P}_\mathscr{B}(\vect{Z})$, Problem \eqref{eq:dual} rewrites as
\begin{equation}
\begin{aligned} 
\argmin_{\vect{Z},\vect{P}\in{\C}} &\quad -\log\det(\hat{\vect{\Sigma}}_{\vec{\mathsf{y}}} +\vect{Z}) - mN\\
\text{subject to } &\quad \eqref{eq:minY1},\,\eqref{eq:minY2}\\
&\quad \vect{P}=\lambda_LI_{mN} + \mathsf{P}_\mathscr{B}(\vect{Z})\\
&\quad \vect{P} \ge 0.
\end{aligned}
\end{equation}
where we have omitted the domain of the objective function $\hat{\vect{\Sigma}}_{\vec{\mathsf{y}}}+\vect{Z}>0$ since it will be checked in the stepsize-choice stage of the algorithm. The augmented Lagrangian for the problem is
\begin{equation*}
\begin{aligned}
\mathcal{L}_\rho(\vect{Z},\vect{P},\vect{M}) = &-\log\det\left(\hat{\vect{\Sigma}}_{\vec{\mathsf{y}}} + \vect{Z}\right) - \left<\vect{M},\,\vect{P}-\lambda_LI_{mN} -\mathsf{P}_\mathscr{B}(\vect{Z})\right>_{\C}\\
&+ \frac{\rho}{2}\,\| \vect{P}-\lambda_LI_{mN} -\mathsf{P}_\mathscr{B}(\vect{Z})\|^2_{\C}
\end{aligned}
\end{equation*}
where $\rho>0$ is the penalty term and $\vect{M}\in{\C}$ is the Lagrange multiplier associated to the equality constraint on $\vect{P}$. Accordingly, the ADMM updates are the following:
\begin{enumerate}
	\item The $\vect{Z}$-minimization step
	\begin{equation}\label{eq:Zupdt}
	\begin{aligned} 
	\vect{Z}^{k+1}=&\argmin_{\vect{Z}\in{\C}}  \quad\mathcal{L}_\rho(\vect{Z},\vect{P}^k,\vect{M}^k)\\
	&\text{subject to } \quad \vect{Z}\in\Z.
	\end{aligned}
	\end{equation}
	\item The $\vect{P}$-minimization step
	\begin{equation}\label{eq:Pupdt}
	\begin{aligned} 
	\vect{P}^{k+1}=&\argmin_{\vect{P}\in{\C}} \quad\mathcal{L}_\rho(\vect{Z}^{k+1},\vect{P},\vect{M}^k)\\
	&\text{subject to } \quad \vect{P}\ge 0.
	\end{aligned}
	\end{equation}
	\item Dual variable update
	\begin{equation}\label{eq:Mupdt}
	\vect{M}^{k+1} = \vect{M}^k - \rho\left(\vect{P}^{k+1}-\lambda_L I_{mN}-\mathsf{P}_\mathscr{B}(\vect{Z}^{k+1})\right).
	\end{equation}
\end{enumerate}
where $\Z:=\{\vect{Z}\in\C:\,\eqref{eq:minY1},\,\eqref{eq:minY2}\}$ and we have considered a constant value of $\rho$ in order simplify the notation. We will discuss later how to update $\rho$ to get a faster convergence. Updates $1)$ and $2)$ are not in an implementable format. The $\vect{Z}$-update step \eqref{eq:Zupdt} is equivalent to the minimization of 
\begin{equation*}
\begin{aligned}
	\I(\vect{Z}) := &-\log\det(\hat{\vect{\Sigma}}_{\vec{\mathsf{y}}}+\vect{Z}) + \frac{\rho}{2}\,\|\mathsf{P}_\mathscr{B}(\vect{Z})\|^2_{\C}\\
&+ \left<\vect{M}^k - \rho\,(\vect{P}^k-\lambda_L I_{mN}),\,\mathsf{P}_\mathscr{B}(\vect{Z})\right>_{\C},
\end{aligned}
\end{equation*}
over the set $\Z$, which has no closed-form solution, as noticed in \cite{SongVan} where the solution is approximated by a projective-gradient step. Following the same lines, the new $\vect{Z}$-update step starts from a known feasible point $\vect{Z}^0=\bar{\vect{Z}}$ and continue the iterations following the update rule
\begin{equation}\label{eq:newZupdt}
\vect{Z}^{k+1} = \mathsf{P}_\Z\left(\vect{Z}^k - t_k\,\nabla\I(\vect{Z}^k)\right)
\end{equation}
where
\[
\nabla\I(\vect{Z}^k) = -(\hat{\vect{\Sigma}}_{\vec{\mathsf{y}}}+\vect{Z}^k)^{-1}+\mathsf{P}_\mathscr{B}(\vect{M}^k) +\rho\,\mathsf{P}_\mathscr{B}\left(\vect{Z}^k-\vect{P}^k + \lambda_L I_{mN}\right)
\]
is the gradient of the cost-function $\I$ computed in $\vect{Z}^k$, $t_k$ is the stepsize founded by the Armijo condition, and $\mathsf{P}_\Z$ is the projection operator onto the constraints space $\Z$.\\
The optimization problem involved in the $\vect{P}$-update step \eqref{eq:Pupdt} is equivalent to minimize the functional
\[
\J(\vect{P}) := \frac{\rho}{2}\,\|\vect{P}\|_{\C}^2 - \left<\vect{P},\,\vect{M}^k+\rho\left(\lambda_L I_{mN}+\mathsf{P}_\mathscr{B}(\vect{Z}^{k+1})\right)\right>_{\C}
\]
over all $\vect{P}\ge 0$. Since $\J$ is a quadratic functional of $\vect{P}$ the  minimization of $\J$ over the whole vector space ${\C}$ admits the closed form solution
\[
\vect{P}_o = \frac{1}{\rho}\,\vect{M}^k+\lambda_L I_{mN} + \mathsf{P}_\mathscr{B}(\vect{Z}^{k+1})
\]
but it is not a positive semidefinite matrix in general. Accordingly, in order to find our solution, we have to find the positive semidefinite block-circulant matrix that better approximates $\vect{P}_o$ in the norm induced by the scalar product on $\C$ (i.e. the Frobenius norm on $\C$). Recall the following well-known result.
\begin{lemma}\label{lem:sdpproj}
	Let $A\in\CC^{n\times n}$ be an Hermitian matrix whose eigenvalue decomposition is given by $A=U^* \Lambda U$, with 
	\[
	U^* U = UU^* = I\qquad\text{and}\qquad \Lambda=\text{diag}\{\lambda_1,\dots,\lambda_n\}.
	\]
	Then, the positive semidefinite matrix that better approximates $A$ in the Frobenius norm is the projection of $A$ onto the cone of positive semidefinite matrices $\mathcal{P}^+$, namely 
	\[
	\mathsf{P}_{\mathcal{P}^+}(A):=\argmin_{X\ge0}\|X-A\|_F = U^*\,\text{diag}\{\gamma^o_1,\dots,\gamma^o_n\}\,U,
	\]
	where
	\[
	\gamma^o_i =
	\left\{
	\begin{aligned}
	&\lambda_{ii},\quad&\text{if}&\quad\lambda_{ii}\ge 0,\\
	&0,\quad&\text{if}&\quad\lambda_{ii}< 0.
	\end{aligned}
	\right.
	\]
\end{lemma}
The following proposition ensures that the projection of a symmetric, block-circulant matrix onto the cone of positive semi-definite matrices is still block-circulant.
\begin{proposition}\label{prop:cpd}
	Let $\vect{C}$ be a symmetric, block-circulant matrix
	\[
	\vect{C} = \vect{F}^*\text{diag}\left\{C(\zeta^0),\,C(\zeta^1),\,\dots,\,C(\zeta^{N-1})\right\}
	\vect{F},
	\]
	and let $C(\zeta^k) = V_k\Lambda_k V_k^*$ with
	\[
		V_k^*V_k = V_kV_k^* = I_m\quad\text{ and }\quad \Lambda_k=\text{diag}\{\lambda_{k1},\dots,\lambda_{km}\},
	\]
	being the eigen-decomposition of the (Hermitian) block $C(\zeta^k)$, for $k=0,\dots,N-1$. Then the eigen-decomposition of $\vect{C}$ can be written as 
	\[
	   \vect{C} = \vect{W}^*\vect{\Lambda}\vect{W},\qquad \vect{W}=\vect{V}^*\vect{F},
	\]
	where $\vect{V} = \text{diag}\{V_0,\dots,V_{N-1}\}$  and $\vect{\Lambda} = \text{diag}\{\Lambda_0,\dots,\Lambda_{N-1}\}$. Then,
	\[
	\mathsf{P}_{\C^+}(\vect{C}) :=\argmin_{\vect{X}\ge0}\|\vect{X}-\vect{C}\|_\C= \vect{W}^*\,\text{diag}\{\Gamma_0,\dots,\Gamma_{N-1}\}\,\vect{W}
	\] 
	where $\Gamma_k = \text{diag}\{\gamma_{k1},\dots,\gamma_{km}\}$ and
	\[
	\gamma_{ki} =
	\left\{
	\begin{aligned}
	&\lambda_{ki},\quad&\text{if}&\quad\lambda_{ki}\ge 0,\\
	&0,\quad&\text{if}&\quad\lambda_{ki}< 0.
	\end{aligned}
	\right.
	\]
	for $k=0,\dots,N-1$.
	\begin{IEEEproof}
	    The result follows from applying Lemma \ref{lem:sdpproj} with $U=\vect{W}$. Of course, 
	    \[
	    \mathsf{P}_{\C^+}(\vect{C}) = \vect{F}^* \text{diag}\{V_0\Gamma_0 V_0^*,\dots,V_{N-1}\Gamma_{N-1} V_{N-1}^*\}
	    \vect{F}
	    \]
	    is a block-circulant matrix because it is block-diagonalized by the Fourier-block matrix.
    \end{IEEEproof}
\end{proposition}

According to Proposition \ref{prop:cpd}, the positive semidefinite block-circulant matrix that better approximates $\vect{P}_o$ in the $\C$ norm is the projection of $\vect{P}_o$ onto the cone of the symmetric, positive semidefinite, block-circulant matrices $\C^+$, that is
\begin{equation}\label{eq:newPupdt}   
\vect{P}^{k+1} = \mathsf{P}_{\C^+}(\vect{P}_o) = \mathsf{P}_{\C^+}\left( \frac{1}{\rho}\,\vect{M}^k+\lambda_L I_{mN} + \mathsf{P}_\mathscr{B}(\vect{Z}^{k+1})\right).
\end{equation}
We conclude that the ADMM algorithm for the estimation of the sparse and the low-rank component of the inverse of the covariance matrix of the reciprocal process consists in the following updates
\begin{equation}\label{eq:ADMMrec}
	\begin{aligned}
	\vect{Z}^{k+1} &= \mathsf{P}_\Z\left[\vect{Z}^k - t_k\,\nabla\I(\vect{Z}^k)\right],\\
    \vect{P}^{k+1} &= \mathsf{P}_{\C^+}\left[ \frac{1}{\rho^k}\,\vect{M}^k+\lambda_L I_{mN} + \mathsf{P}_\mathscr{B}(\vect{Z}^{k+1})\right],\\
    \vect{M}^{k+1} &= \vect{M}^k - \rho^k\left[\vect{P}^{k+1}-\lambda_L I_{mN}-\mathsf{P}_\mathscr{B}(\vect{Z}^{k+1})\right].
	\end{aligned}
\end{equation}
A typical update for $\rho$ is $\rho^{k+1}=\alpha\rho^k$, with $\alpha>1$ being a certain growth coefficient that needs to be properly tuned.
Notice that the matrices involved in (\ref{eq:ADMMrec}) are all symmetric and block-circulant. Accordingly, as explained in Section \ref{sec:rp}, the introduction of the reciprocal approximation allows to obtain a robust identification procedure even in the case when $n$ is large. Indeed, relations \eqref{eq:covrec} and \eqref{eq:invcovrec}, allow to compute inverse matrices and eigenvalues in a robust way. Moreover, it is worth noting from \eqref{eq:invcovrec}, that the dimensions of the matrices whose eigenvalues must be computed in the optimization procedure, depend only on $m$, hence the identification algorithm we are proposing scales with respect to $n$ gaining robustness in the results even if the order of the AR process is large.\\
Following \cite{BoydADMM}, the basic stopping criterium for the algorithm is based on the primal and dual residuals of the optimality conditions that respectively measure the satisfaction of the inequality constraint $\vect{P}\ge0$ and the distance between two successive iterates of the variable $\vect{P}$. More precisely, the primal residual at iteration $k+1$ is defined as
\[
r^{k+1} := \vect{P}^{k+1}-\lambda_L I_{mN} - \mathsf{P}_\mathscr{B}(\vect{Z}^{k+1}),
\]
while the dual residual turns out to be
\[
s^{k+1} := \mathsf{P}_{\mathscr{B}^c}(\vect{M}^k) - \rho^k\,\left[\vect{P}^{k+1} - \mathsf{P}_\mathscr{B}(\vect{P}^k) \right].
\]
It is reasonable that the primal and dual residual must be small, that is
\[
\|r^k\|_{\C}\le\epsilon^\text{p}\qquad\text{and}\qquad \|s^k\|_{\C}\le\epsilon^\text{d},
\]
where $\epsilon^\text{p}>0$ and $\epsilon^\text{d}>0$ are feasibility tolerances for the primal and dual feasibility conditions. The latter are defined as
\begin{align*}
\epsilon^\text{p} &:= mN\,\epsilon^\text{abs} + \epsilon^\text{rel}\max\left\{\lambda_L\sqrt{mN},\,\|\vect{Z}^k\|_{\C},\,\|\vect{P}^k\|_{\C}\right\},\\
\epsilon^\text{d} &:= mN\,\epsilon^\text{abs} + \epsilon^\text{rel}\,\|\vect{M}^k\|_{\C}.
\end{align*}
Here, $\epsilon^\text{abs}$ and $\epsilon^\text{rel}$ are predefined absolute and relative tolerances for the problem. Accordingly, the algorithm converges if all the conditions
\begin{equation}\label{eq:stopcrit}
\|r^k\|_{\C}\le\epsilon^\text{p},\qquad \|s^k\|_{\C}\le\epsilon^\text{d},\qquad
\rho^k=\rho_\text{max}
\end{equation}
hold true, where $\rho_\text{max}>0$ is the maximum value allowed for the penalty parameter $\rho^k$, selected by the user.

\section{NUMERICAL EXAMPLES}\label{sec:ne}
In this section we compare the performances of our method to which we will refer to as \emph{approximated algorithm} with the one proposed in \cite{ZorziADMM} for the solution of Problem \eqref{eq:SpLcl}, which will be referred to as \emph{exact algorithm}. In particular we will show how the two algorithms behave considering both the case in which the observed process has low dimension and the case in which we have an high dimensional observed process.
\subsubsection*{Low-dimensional case} 
\begin{figure}[h!]
	\centering
	\includegraphics[scale=0.34]{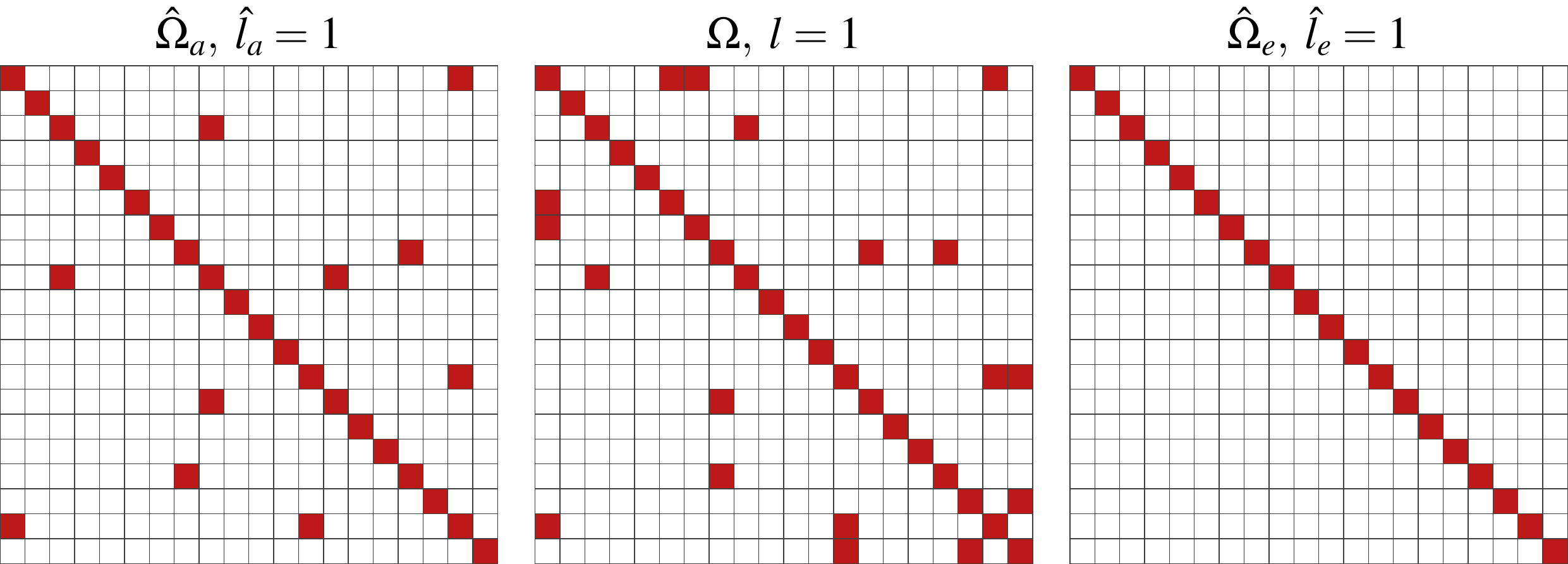}
	\caption{Sparsity pattern estimated by the approximated algorithm with $\alpha=1.007,\,\lambda_S=95,\,\lambda_L=5.4$ (left), true sparsity pattern (center), sparsity pattern estimated by the exact algorithm with $\alpha=1.002,\,\gamma_S=2.6,\,\gamma_L=2.95$ (right). The red squares indicate the conditional dependent pairs while the white squares indicates the conditional independent pairs. $\hat l_a$, $l$ and $\hat l_e$ denote the number of latent variables.}
	\label{fig:truemod}
\end{figure}

\begin{figure}[h!]
	\centering
	\includegraphics[scale=0.195]{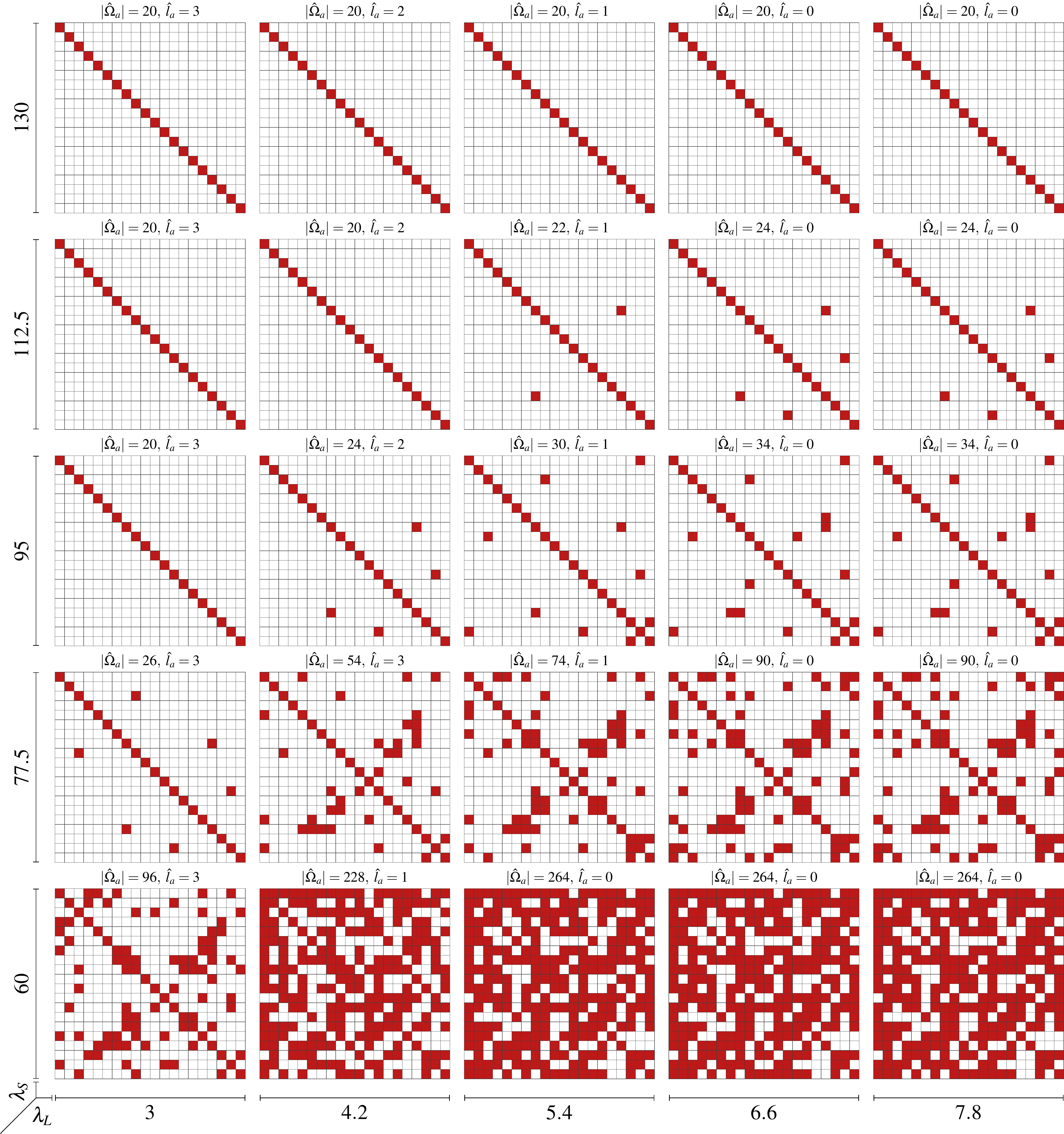}
	\caption{Supports and ranks estimated by the approximated algorithm for $\lambda_S\in[60,130]$ and $\lambda_L\in[3,7.8]$. The growth coefficient is set $\alpha=1.007$.}
	\label{fig:allmod}
\end{figure}

\begin{figure}[h!]
	\centering
	\includegraphics[scale=0.53]{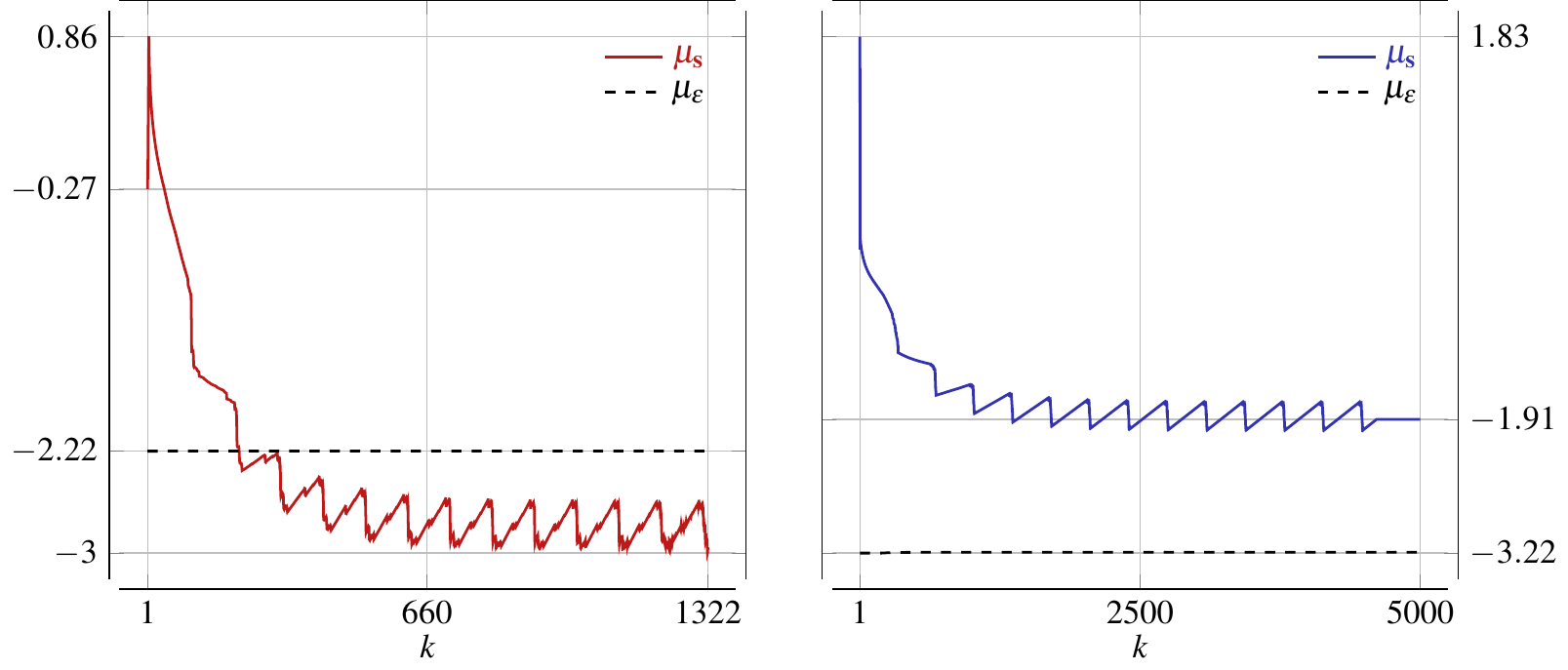}
	\caption{Logarithm of the average dual residual for the approximated method (left) and for the exact method (right). The dashed lines correspond to the logarithm of the associated average feasibility tolerances.}
	\label{fig:resid}
\end{figure}

\begin{figure}[h!]
	\centering
	\includegraphics[scale=0.55]{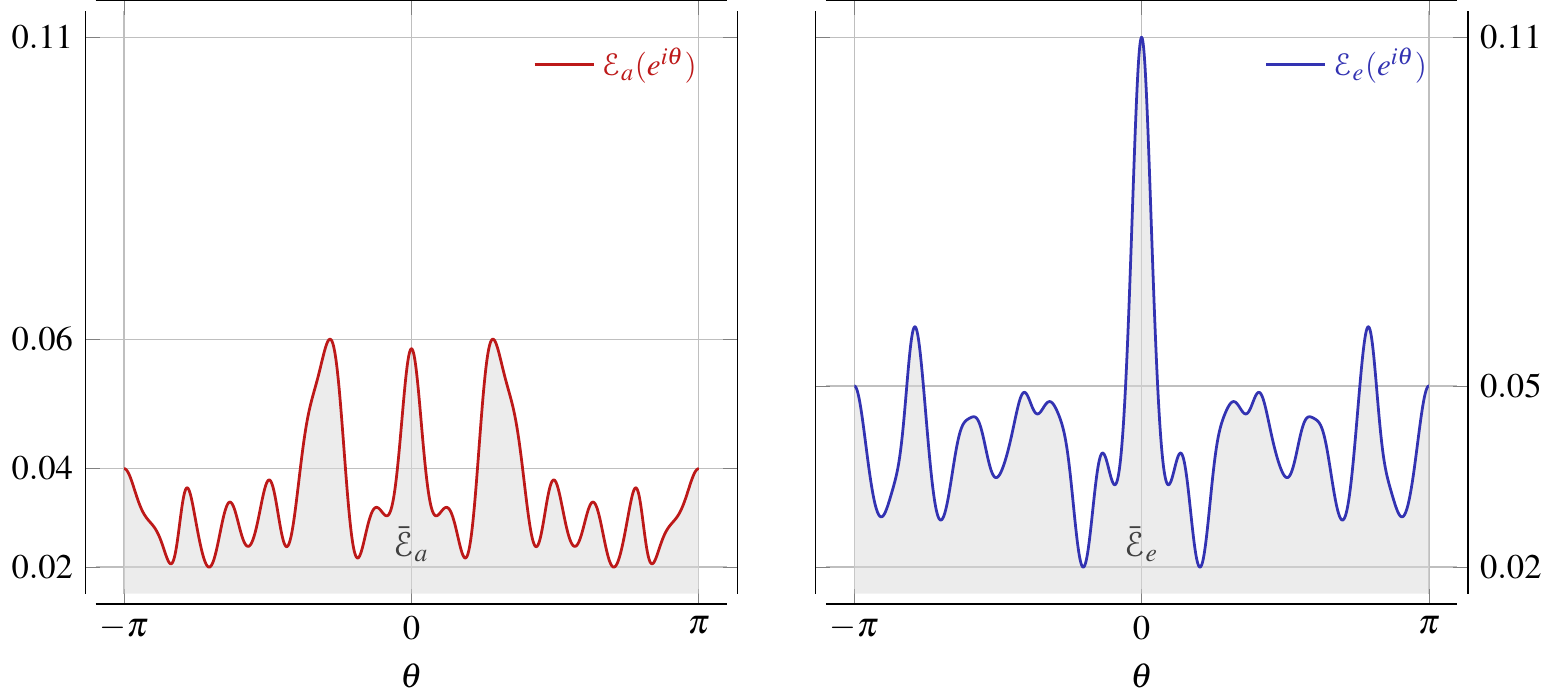}
	\caption{Relative errors in the estimated spectra: approximated algorithm (left), exact algorithm (right).}
	\label{fig:DRpsd}
\end{figure}

Synthetic data ere generated from the AR latent-variable model of order $n=8$,
\begin{equation}\label{eq:ARsim}
   \vec{\mathsf{y}}(t)=\sum_{k=1}^n A_k\,\vec{\mathsf{y}}(t-k)+\vec{\eta}(t),
\end{equation}
with $m=20$ observed variables and $l=1$ latent variables. Here, $\vec{\eta}(t)$ is white Gaussian noise with variance  $\EE[\vec{\eta}(t)^\top\vec{\eta}(t)]=21.14$ and $T=1000$ samples have been used to compute the estimated covariance lags $\hat{R}_k$, $k=0,\dots,n$. Figure \ref{fig:truemod} (center) reports the sparsity pattern of the underlying model, randomly generated so that the non-zero elements represents the $5\%$ of the total elements.
For the approximated algorithm we have considered $N=30$ samples of the spectrum. In both the ADMM implementations we have set $\epsilon^\text{abs}=10^{-5}$ and $\epsilon^\text{rel}=10^{-4}$ while $\rho_\text{max}=10^4$. In order to tune the update of the penalty term $\rho$ in the ADMM, we have ran both the algorithms for different values of the growth coefficient $\alpha\in[1.001,\,1.1]$. More precisely, for each value of $\alpha$, a $5\times 5$ grid of candidate estimated models has been produced, corresponding to five linearly spaced values of the regularization parameters $\lambda_S\in[60,130]$ and $\lambda_L\in[3,7.8]$ for the approximated algorithm, and five linearly spaced values of $\gamma_S\in[1.42,2.6]$ and $\gamma_L\in[2.425,2.95]$ for the exact algorithm. The values of the regularization parameters that identify the grids have been selected so that the estimated models capture a range of features as complete as possible: from a very sparse model with a relatively high rank, to a quasi-full model with the lowest rank possible. Figure \ref{fig:allmod} shows the supports and the ranks estimated by the approximated algorithm corresponding to the different values of $\lambda_S$ and $\lambda_L$. For both methods the value of $\alpha$ that gives the better performances, i.e. that guarantees the minimum gap between $\epsilon^\text{p}\, /\,  \epsilon^\text{d}$ and the primal/dual residual at the final iteration, respectively, has been selected. Accordingly, we have chosen $\alpha=1.007$ for the approximated algorithm while $\alpha=1.002$ has been chosen for the exact algorithm.
Let $\vect{s}(h)$ and $\vect{\epsilon}(h)$ denote the vectors containing the dual residual and its feasibility tolerance for the model $h=1,\dots,25$ respectively. Figure \ref{fig:resid} displays the (logarithm of the) averages
\[
 \mu_\vect{s} = \frac{1}{25}\sum_{h=1}^{25} \vect{s}(h),
 \qquad
 \mu_{\vect{\epsilon}} = \frac{1}{25}\sum_{h=1}^{25} \vect{\epsilon}(h),
\]
obtained by our method with $\alpha=1.007$ (left) and by the exact method for $\alpha=1.002$ (right). For both algorithms, the primal residual always satisfies the condition in the stopping criterium \eqref{eq:stopcrit} therefore there is no need to displaying it.
We observe that the exact algorithm does not converge for any value of $\alpha$ we have considered. Indeed, the plot in Figure \ref{fig:resid} (right) clearly shows that the mean dual-residual $\mu_\vect{s}$ stays significantly above the threshold $\mu_{\vect{\epsilon}}$. The optimal values of the regularization parameters have then been selected by cross-validation, using a test data set of $500$ samples. Figure \ref{fig:truemod} compares the optimal sparsity pattern provided by the approximated algorithm $\hat{\Omega}_a$ (left), corresponding to $\lambda_S=95$ and $\lambda_L=5.4$, and the optimal sparsity pattern estimated by the exact algorithm $\hat{\Omega}_e$ (right) corresponding to $\gamma_S=2.6$ and $\gamma_L=2.95$, together with the estimates of the number of latent variables, $\hat{l}_a$ and $\hat{l}_e$, respectively. Notice that both algorithms estimates the correct number of latent variables but only the approximated one produces an estimate of the sparsity pattern comparable with the true one. Let $\hat{\Phi}_e$ and $\hat{\Phi}_a$ be the estimates of the spectra of the true observed process $\Phi_{\vec{\mathsf{y}}}$ obtained by the solutions of problems \eqref{eq:SpLcl} and \eqref{eq:primal}, respectively. According to Figure \ref{fig:recAR}, $\hat{\Phi}_a$ is the extension over the whole interval $[-\pi,\pi]$ of the symbol of the estimated covariance matrix of the reciprocal process $\vect{y}$ approximating $\vec{\mathsf{y}}$. The squared-estimation errors for the two algorithms are depicted in Figure \ref{fig:DRpsd}; the corresponding mean values over $[-\pi,\pi]$ are 
\begin{align*}
\mathcal{E}_a&:=\dfrac{\|\Phi_{\vec{\mathsf{y}}}-\hat{\Phi}_a\|_F^2}{\|\Phi_{\vec{\mathsf{y}}}\|_F^2},&\quad\bar{\mathcal{E}}_a&:=\int\,\mathcal{E}_a(e^{i\theta})=0.0358,\\
\mathcal{E}_e&:=\dfrac{\|\Phi_{\vec{\mathsf{y}}}-\hat{\Phi}_e\|_F^2}{\|\Phi_{\vec{\mathsf{y}}}\|_F^2},&\quad\bar{\mathcal{E}}_e&:=\int\,\mathcal{E}_e(e^{i\theta})=0.0443.
\end{align*}
The approximated algorithm performs better both in terms of the mean value and in terms of the height of the peaks of the relative error.   
\subsubsection*{High-dimensional case} 
\begin{figure}[h!]
	\centering
	\includegraphics[scale=0.66]{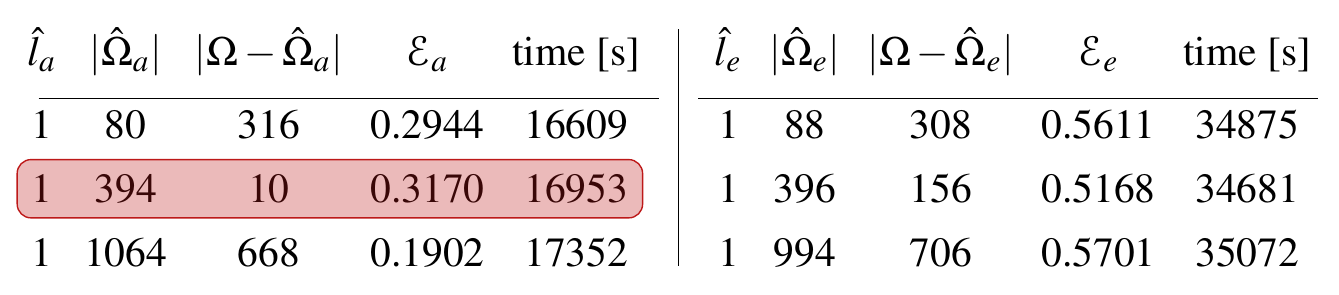}
	\caption{Summary of the performances of the two algorithm for $\lambda_S=100,\,146.25,\,350$ and $\gamma_S=0.7,\,0.826,\,1.7$. The values of the low-rank regularization parameters are $\lambda_L=8.6875$ for the approximated algorithm (left) and $\gamma_L=2.3$ for the exact algorithm (right). These results have been obtained on a 2014 1.4GHz MacBook Air.}\label{fig:tabHD}
\end{figure}

We consider now an AR latent-variable model as in \eqref{eq:ARsim} where we have $m=80$ observed variables and $l=1$ latent variable, $n=24$ and the variance of the noise is $\EE[\vec{\eta}(t)^\top\vec{\eta}(t)]=87.1$. The number of samples used to estimate the covariance lags $R_k$ is $T=15000$. The number of conditionally dependent pairs in the true model is $158$ so that the cardinality of the true support is $|\Omega|=396$. Table \ref{fig:tabHD} compares the performances of our approximated algorithm with the exact algorithm proposed in \cite{ZorziADMM} for different values of the sparsity regularization parameters $\lambda_S$ and $\gamma_S$, that have been chosen in order to have approximatively the same variety on the results. The notation $|\Omega -\hat{\Omega}|$ indicates the error on the sparsity pattern in terms of number of misclassified entries.
Both algorithms estimate the correct number of latent variables, but the approximated algorithm gives a result very close to the true one (highlighted in red in Figure \ref{fig:tabHD}) while for the exact algorithm, even if the cardinality of the true support has been correctly estimated, the error in the reconstruction of the sparsity pattern is quite high. This is due to the fact that the higher is the order of the process $n$, the less accurate is the computation of eigenvalues and inverse matrices by the exact algorithm. Figure \ref{fig:tabHD} shows that such an issue is avoided in the approximated version, thanks to the availability of closed-form formulas for the computation of the eigenvalues of block-circulant matrices. Moreover, we see that the run time of the exact algorithm is about twice the run time of the approximated one. This confirm the fact that the approximated algorithm scales with the order $n$ of the AR process we are approximating as suggested in Section \ref{sec:id}. This kind of scenario agrees with what we have discussed in Section \ref{sec:id}: high-order AR process are quite challenging instances for the exact procedure proposed in \cite{ZorziADMM}; in this cases, the reciprocal approximation leads to remarkable benefits in the performances of the identification procedure.

\section{CONCLUSIONS}\label{sec:conc}
In this paper an identification paradigm for latent-variable graphical models associated to reciprocal processes has been presented. It has been shown that the proposed paradigm is theoretically strongly sustained, being an approximation of the corresponding problem for AR processes both in a maximum likelihood and in a maximum entropy sense. The performances of the proposed method have been compared with the approach proposed in \cite{ZorziADMM} where no approximation is introduced. The numerical examples have shown that for high-order AR processes reciprocal approximation gives substantial improvements in terms of robustness and scalability of the identification procedure.

\bibliography{biblio}
\bibliographystyle{ieeetr}

\end{document}